\documentclass[12pt,reqno]{amsart}

\usepackage{amscd}
\usepackage{epsfig}
\usepackage{float}
\usepackage{url}

\newcommand{\N}{{\mathbb N}}
\newcommand{\Z}{{\mathbb Z}}
\newcommand{\Q}{{\mathbb Q}}
\newcommand{\C}{{\mathbb C}}
\newcommand{\R}{{\mathbb R}}
\renewcommand{\P}{{\mathbb P}}

\newcommand{\FF}{{\mathcal F}}

\newcommand{\ddd}{{\rm d}}
\newcommand{\paa}{{\partial}}

\begin{document}

\theoremstyle{plain}
\newtheorem{lemma}{Lemma}[section]
\newtheorem{theorem}[lemma]{Theorem}
\newtheorem{proposition}[lemma]{Proposition}
\newtheorem{corollary}[lemma]{Corollary}
\newtheorem{conjecture}[lemma]{Conjecture}
\newtheorem{conjectures}[lemma]{Conjectures}

\theoremstyle{definition}
\newtheorem{definition}[lemma]{Definition}
\newtheorem{withouttitle}[lemma]{}
\newtheorem{remark}[lemma]{Remark}
\newtheorem{remarks}[lemma]{Remarks}
\newtheorem{example}[lemma]{Example}
\newtheorem{examples}[lemma]{Examples}

\title[Bernoulli moments of spectral numbers]
{Bernoulli moments of spectral numbers and Hodge numbers} 

\author{Thomas Br\'elivet\and Claus Hertling}

\address{Thomas Br\'elivet\\Departamento de \'Algebra, Geometr\'ia 
y Topolog\'ia, Facultad de ciencias, Universidad de Valladolid,
47005 Valladolid, Espa\~na}

\email{brelivet\char64 agt.uva.es}

\address{Claus Hertling\\
Institut \'Elie Cartan (Math\'ematiques), Universit\'e Henri Poincar\'e,
B.P. 239, 54506 Vandoeuvre-l\`es-Nancy, France}

\email{hertling\char64 iecn.u-nancy.fr}




\begin{abstract}
The distribution of the spectral numbers of an isolated hypersurface
singularity is studied in terms of the Bernoulli moments. These
are certain rational linear combinations
of the higher moments of the spectral numbers. They are related to the
generalized Bernoulli polynomials. We conjecture that their signs 
are alternating and prove this in many cases. One motivation for
the Bernoulli moments comes from the comparison with compact complex
manifolds.
\end{abstract}

\maketitle

\tableofcontents

\setcounter{section}{0}

\section{Conjectures and results}\label{c1}
\setcounter{equation}{0}

\noindent
An isolated hypersurface singularity $f:(\C^{n+1},0)\to (\C,0)$
with Milnor number $\mu$ 
comes equipped with its spectral numbers, a tuple of $\mu$
rational numbers $\alpha_1,...,\alpha_\mu$ which satisfy
\begin{eqnarray}\label{1.1}
-1<\alpha_1\leq ... \leq \alpha_\mu <n\quad \hbox{ and } 
\quad \alpha_i+\alpha_{\mu+1-i} = n-1.
\end{eqnarray}
They come from the Hodge filtration on the middle cohomology of
the Milnor fiber and the semisimple part of the monodromy, acting
on it \cite{St}\cite{AGV}.

We are interested in their distribution. We consider the numbers
\begin{eqnarray}\label{1.2}
V_{2k}^{sing}(f):= \sum_{i=1}^\mu \left(\alpha_i-\frac{n-1}{2}\right)^{2k}
\quad\hbox { for } k\geq 0.
\end{eqnarray}
One should divide them by $\mu=V_0^{sing}(f)$ to get the normalized moments, 
but we prefer not to do it.
So we call $V_2^{sing}(f)$ the {\it variance} 
of $f$ and the $V_{2k}^{sing}(f)$ the {\it higher moments}.
In \cite{He1}\cite{He2} C. Hertling formulated the conjecture that 
any isolated hypersurface singularity satisfies
\begin{eqnarray}\label{1.3}
V_2^{sing}(f)\leq V_0^{sing}(f)\cdot \frac{\alpha_\mu-\alpha_1}{12}.
\end{eqnarray}
It was proved by M. Saito for irreducible curve singularities \cite{SM2},
by T. Br\'elivet for curve singularities with nondegenerate Newton
boundary \cite{Br1}, and recently by T. Br\'elivet for all curve
singularities \cite{Br3}. A. Dimca formulated a dual conjecture
with $\geq$ instead of $\leq $ for polynomials with isolated 
singularities \cite{Di}.

In this paper, the conjecture \eqref{1.3} will be extended to a series
of inequalities for certain linear combinations of the 
higher moments, which will be called {\it Bernoulli moments}. 
Before explaining this, we regard a related situation,
where these linear combinations will also be interesting.

If $X$ is a compact complex manifold of dimension $n$, one often considers
\cite{Hi}
\begin{eqnarray}\label{1.4}
h^{p,q}&=& \dim H^q(X, \Omega^p)\qquad \hbox{ \ and }\\
\chi_p &=& (-1)^p\chi(\Omega^p) = \sum_q(-1)^{p+q}h^{p,q} \label{1.5}
\end{eqnarray}
We define
\begin{eqnarray}\label{1.6}
V_{2k}^{mfd}(X):=\sum_p \chi_p(p-\frac{n}{2})^{2k} \quad \hbox{ for }k\geq 0.
\end{eqnarray}
At the end of this chapter and in the last chapter, 
we will consider this situation.
Here $V_0^{mfd}(X)=\int_X c_n$ could be 0; then we cannot normalize
the moments $V_{2k}^{mfd}(X)$.

Now let 
\begin{eqnarray}\label{1.7}
V = \sum_{k=0}^\infty V_{2k}\frac{1}{(2k)!}t^{2k}
\end{eqnarray}
be a formal power series in $t^2$ with variables $V_0,V_2,V_4,...$,
and let $\nu$ be another variable.
The Bernoulli moments $\Gamma^{Ber}_{2k}(V,\nu)\in \sum_{l=0}^k\Q[\nu]V_{2l}$
are defined by
\begin{eqnarray}\label{1.8}
\Gamma^{Ber}(V,\nu)=\sum_{k=0}^\infty 
\Gamma_{2k}^{Ber}\frac{1}{(2k)!}t^{2k}
= V\cdot \exp\left( \nu\cdot \log\frac{t/2}{\sinh (t/2)}\right).&&
\end{eqnarray}

The first four of them are 
\begin{eqnarray}\label{1.9}
\Gamma_0^{Ber}(V,\nu)&=&V_0,\\
\Gamma_2^{Ber}(V,\nu)&=&V_2-V_0\cdot(\frac{1}{12}\nu),\label{1.10}\\
\Gamma_4^{Ber}(V,\nu)&=&V_4 - V_2\cdot(\frac{1}{2}\nu) + 
V_0\cdot(\frac{1}{120}\nu +\frac{1}{48}\nu^2), \label{1.11}\\
\Gamma_6^{Ber}(V,\nu)&=&V_6 - V_4\cdot(\frac{5}{4}\nu) 
+ V_2\cdot(\frac{1}{8}\nu + \frac{5}{16}\nu^2) \nonumber\\
&&- V_0\cdot(\frac{1}{252}\nu + \frac{1}{96}\nu^2 + \frac{5}{576}\nu^3).
\label{1.12}
\end{eqnarray}

The Bernoulli moments are closely related to the generalized Bernoulli
polynomials. This will be discussed after theorem \ref{t1.3}.
A relation with the Bernoulli numbers $B_n$ is given by
\begin{eqnarray}\label{1.13}
\log\frac{t/2}{\sinh (t/2)}
= \sum_{k=1}^\infty \frac{-1}{2k}B_{2k} \frac{1}{(2k)!}t^{2k}.
\end{eqnarray}
The Bernoulli numbers $B_{2k}$ for $k\geq 1$ satisfy 
$B_{2k}\in (-1)^{k-1}\Q_{>0}$. Therefore the coefficient of $V_{2j}$
in $\Gamma^{Ber}_{2k}$ is a polynomial in $\nu$ 
with all coefficients having the sign $(-1)^{k-j}$. A more precise
discussion in chapter \ref{c2} shows the following elementary lemma.

\begin{lemma}\label{t1.1}
Consider $V=\sum_{k=0}^{\infty}V_{2k}\frac{1}{(2k)!}t^{2k}\in \R[[t]]$ 
with $V_0>0$. Fix $k_0\in \N\cup\{\infty\}$
(with $\N:=\{0,1,2,...\}$ in this paper).
\begin{list}{}{}
\item[a)] If $k_0\in \N$, there exists a number $\nu\in \R$ such that
\begin{eqnarray}\label{1.14}
\forall\ k\in \N \hbox{ with }k \leq k_0\qquad 
(-1)^k\Gamma_{2k}^{Ber}(V,\nu)\geq 0.
\end{eqnarray}
\item[b)] If a number $\nu\in\R$ satisfies \eqref{1.14} for 
$k_0\in \N\cup\{\infty\}$ then
also any number $\nu'\in \R$ with $\nu'>\nu$ satisfies \eqref{1.14}.
\end{list}
\end{lemma}

In view of this lemma, the first of the following two conjectures
is weaker than the second. These conjectures are at the center  
of this paper.

\begin{conjectures}\label{t1.2}
Let $f:(\C^{n+1},0)\to (\C,0)$ be an isolated hypersurface singularity.
\begin{list}{}{}
\item[(W)] (Weak form) Then for all $k\in \N$
\begin{eqnarray}\label{1.15}
(-1)^k\Gamma_{2k}^{Ber}(V^{sing}(f),{n+1})>0.
\end{eqnarray}
\item[(S)] (Strong form) Then for all $k\in \N$
\begin{eqnarray}\label{1.16}
(-1)^k\Gamma_{2k}^{Ber}(V^{sing}(f),{\alpha_\mu-\alpha_1})\geq 0.
\end{eqnarray}
\end{list}
\end{conjectures}

The case $k=1$ of the conjecture (S) is \eqref{1.3}.
Our evidence for the conjectures is collected in the following theorem.

\begin{theorem}\label{t1.3}
a) The conjecture \ref{t1.2} (S) is true for all quasihomogeneous 
singularities.

b) The conjecture \ref{t1.2} (S) is true for all hyperbolic singularities
$T_{pqr}$.

c) The conjecture \ref{t1.2} (W) is true for all irreducible curve 
singularities.

d) \cite{Br3} The conjecture \eqref{1.3} is true for all curve singularities.

e) If the conjecture (S) [respectively the conjecture (W)] is true for
two singularities $f(x)$ and $g(y)$ then it 
is also true for the sum $f(x)+g(y)$.

f) For any singularity $f$ and any $\nu\in \R_{>0}$ a bound $k_0\geq 0$ exists
such that for all $k\in \N$ with $k\geq k_0$
\begin{eqnarray}\label{1.17}
(-1)^k\Gamma_{2k}^{Ber}(V^{sing}(f),\nu)>0.
\end{eqnarray}
\end{theorem}

Part a) and b) will be proved in chapter \ref{c5}. There we will give
also precise formulas. They even suggest to consider the 
$\Gamma_{2k}^{Ber}(V^{sing}(f),\nu)$ for $\nu={\alpha_\mu-\alpha_1}$
and $\nu={n+1}$ themselves as generalisations of the Bernoulli numbers,
two series for each singularity.

Part c) will be proved in chapter \ref{c6}. Part e) is an easy consequence
of the Thom-Sebastiani formula for spectral numbers; it will be discussed
in chapter \ref{c2}, remark \ref{t2.5} b). 
The proof of part f) will be given after theorem \ref{t1.4}.

The generating function $V^{sing}(f)$ of the higher moments $V_{2k}^{sing}(f)$
takes a very special form,
\begin{eqnarray}\label{1.18}
V^{sing}(f)=\sum_{i=1}^\mu \cosh\left(t(\alpha_i-\frac{n-1}{2})\right)
=\sum_{i=1}^\mu e^{t(\alpha_i-\frac{n-1}{2})}.
\end{eqnarray}
The second equality follows from the symmetry in \eqref{1.1}.
This formula and \eqref{1.8} motivate the following definition.

The polynomials $A_k(x,\nu)\in \Q[x,\nu]$ for $k\in \N$ are defined by
\begin{eqnarray}\label{1.19}
e^{xt}\cdot \exp\left( \nu\cdot \log\frac{t/2}{\sinh (t/2)}\right)
= \sum_{k=0}^\infty A_k(x,\nu)\frac{1}{k!}t^{k}.
\end{eqnarray}
Up to a shift in $x$ they are the generalized Bernoulli polynomials,
which were defined by N\"orlund \cite{No3}\cite{No4}\cite{No1}\cite{No2}. 
$A_{k}(x,\nu)$ is a polynomial of degree
$k$ in $x$ and of degree $\left[\frac{k}{2}\right]$ in $\nu$.
The classical Bernoulli polynomials are the polynomials
$B_k(x)= A_k(x-\frac{1}{2},1)$. The Bernoulli numbers are 
$B_k:=B_k(0)$. The polynomials in $x$ for fixed $\nu\in \N$ and especially
for $\nu=1$ have been studied since very long time. 
We review some properties of the $A_k(x,\nu)$ in chapter \ref{c3}. 

\eqref{1.8}, \eqref{1.18} and \eqref{1.19} together show
\begin{eqnarray}\label{1.20}
(-1)^k \Gamma^{Ber}_{2k}(V^{sing}(f),\nu)
= \sum_{j=1}^\mu (-1)^kA_{2k}(\alpha_j-\frac{n-1}{2},\nu).
\end{eqnarray}
This justifies the name Bernoulli moments.
One crucial property of the polynomials $A_k(x,\nu)$ is the following.

\begin{theorem}\cite{No4}\label{t1.4}
On any compact intervall $I\subset \R$ and for any $\nu\in \R-\Z_{\leq 0}$,
the sequence of polynomials 
\begin{eqnarray}\label{1.21}
(-1)^k A_{2k}(x,\nu)\cdot \frac{(2\pi)^{2k}\cdot \Gamma(\nu)}
{2\cdot (2k)!\cdot (2k)^{\nu-1}}
\end{eqnarray}
tends uniformly to $\cos(2\pi x)$ as $k\to \infty$
(here $\Gamma$ is the gamma function).
\end{theorem}

Therefore for any $\nu\in \R-\Z_{\leq 0}$ the sequence of numbers
\begin{eqnarray}\label{1.22}
(-1)^k \Gamma^{Ber}_{2k}(V^{sing}(f),\nu)\cdot 
\frac{(2\pi)^{2k}\cdot \Gamma(\nu)} {2\cdot (2k)!\cdot (2k)^{\nu-1}}
\end{eqnarray}
tends with $k\to \infty$ to 
\begin{eqnarray}\nonumber
\sum_{j=1}^\mu \cos(2\pi (\alpha_j-\frac{n-1}{2}))
&=& (-1)^{n-1}\cdot \sum_{j=1}^\mu e^{2\pi i\alpha_j}\nonumber\\
&=& (-1)^{n-1} \hbox{trace (monodromy)}\nonumber\\
&=& 1.\label{1.23}
\end{eqnarray}
The last equality is a result of A'Campo \cite{AC1}\cite{AC2}.
For $\nu>0$ the factor 
on the right hand side
in \eqref{1.22} is positive. This shows part f) of theorem \ref{t1.3}.

\begin{remarks}\label{t1.5}
i) If the conjecture (W) is true for a singularity $f$, one can define
a sequence of numbers $\nu_k>0$ for $k\in \N$ such that $\nu_k$ is minimal
with the property
\begin{eqnarray}\label{1.24}
\forall\ k'\geq k \ \forall\ \nu'> \nu_k\ 
(-1)^k \Gamma^{Ber}_{2k'}(V^{sing}(f),\nu')>0.
\end{eqnarray}
In view of part f) of theorem \ref{t1.3}, this decreasing sequence tends
to 0. One could ask how fast.

ii) The conjectures \ref{t1.2} and theorem \ref{t1.4} together
predict the sign of the trace of the monodromy. It is an integer.
By A'Campo's result it is the smallest integer with the correct sign.
In view of this, it seems that 
the values $(-1)^k \Gamma^{Ber}_{2k}(V^{sing}(f),\nu)$ are 
``rather small'', up to the factor in \eqref{1.22}.
\end{remarks}

The conjecture \ref{t1.2} (W) is connected
with some work of K. Saito \cite{SK2} on the spectral numbers. 
He defined the associated distribution
\begin{eqnarray}\label{1.25}
\Delta(f)(s):=\sum_{j=1}^\mu \delta(s-\alpha_j+\frac{n-1}{2}),
\end{eqnarray}
where $\delta(s)$ is the delta function. Because of \eqref{1.18}, 
its Fourier transform is just
$V^{sing}(f)(2\pi i t)$. K. Saito proposed to compare $\Delta(f)(s)$ 
with the continuous distribution
\begin{eqnarray}\label{1.26}
\Delta^{(n+1)}(s) := (\Delta^{(1)} \ast ... \ast\Delta^{(1)})(s)
\end{eqnarray}
(the convolution n+1 times), where
\begin{eqnarray}
\Delta^{(1)}(s):= \left\{ 
\begin{array}{rrr}1&\hbox{ if }&s\in [-\frac{1}{2},\frac{1}{2}],\\
0&\hbox{ if }&s\notin [-\frac{1}{2},\frac{1}{2}] \end{array}\right.
\label{1.27}
\end{eqnarray}
He proved that $\Delta^{(n+1)}(s)$ 
is the limit distribution for quasihomogeneous singularities
if the weights tend to zero and for irreducible curve singularities with
$g$ Puiseux pairs if the last denominator tends to infinity.

The Fourier transform of $\Delta^{(n+1)}(s)$ is 
\begin{eqnarray}\label{1.28}
\left(\frac{\sin (\pi t)}{\pi t}\right)^{n+1}.
\end{eqnarray}
Therefore 
\begin{eqnarray}\label{1.29}
\Gamma^{Ber}(V^{sing}(f),n+1)(2\pi it)
= V^{sing}(f)(2\pi it)/ \left(\frac{\sin (\pi t)}{\pi t}\right)^{n+1}
\end{eqnarray}
is the quotient of the Fourier transforms of the actual distribution
$\Delta(f)(s)$ of spectral numbers and the continuous distribution
$\Delta^{(n+1)}(s)$. The conjecture \ref{t1.2} (W) simply predicts that
all its coefficients are positive. In this sense it confirms  a feeling of  
K. Saito \cite[p 202, (2.5) ii)]{SK2} that the limit distribution
$\Delta^{(n+1)}(s)$ should not only be a limit, but also a bound
for the actual distributions. 

But it is difficult to derive from
this conjecture on the Fourier transforms concrete information
on the distribution $\Delta(f)(s)$. It does not imply a conjecture
of K. Saito \cite[p 203]{SK2} (and Durfee in the case $n=2$)
on the number of spectral numbers in $]-1,0]$. We discuss this 
in chapter \ref{c4}.

\bigskip
We presented ample evidence that the Bernoulli moments are natural
objects. A characterisation in corollary \ref{t2.3} and the explicit
formulas in chapter \ref{c5} will even strengthen this.

But we found the Bernoulli moments in a different way, 
by looking at the moments $V_{2k}^{mfd}(X)$ of 
compact complex manifolds. In chapter \ref{c7} the following 
results will be proved, using the Hirzebruch\--Riemann\--Roch theorem.

\begin{theorem}\label{t1.6}
a) There exist polynomials $q_{kj}(\nu,y_1,...,y_j)\in \Q[\nu,y_1,...,y_j]$
for $k\geq 1$ and $0\leq j\leq 2k-1$ with the following properties.
They are quasihomogeneous of degree $j$ with respect to the weights
$i$ of $y_i$. They satisfy $\deg_\nu q_{k0}=k$ and 
$\deg_\nu q_{kj}\leq k-1-\left[\frac{j}{2}\right]$ for $j\geq 1$.
For any compact complex manifold $X$ of any dimension $n$,
\begin{eqnarray}\label{1.30}
V_{2k}^{mfd}(X) = \sum_{j=0}^{\min(2k-1,n)}
\int_X q_{kj}({n},c_1,...,c_j)
\cdot c_{n-j}
\end{eqnarray}
if $k\geq 1$ and $V_{0}^{mfd}(X)=\int_Xc_n$.

b) The Bernoulli moments of $V^{mfd}(X)$ with $\nu=n$ are 
\begin{eqnarray}\label{1.31}
\Gamma_{2k}^{Ber}(V^{mfd}(X),n) = 
\sum_{j=0}^{\min(2k-1,n)} \int_X q_{kj}(0,c_1,...,c_j)\cdot c_{n-j}
\end{eqnarray}
if $k\geq 1$ and $\Gamma_{0}^{Ber}(V^{mfd}(X),n)=\int_Xc_n$.
\end{theorem}

The formulas for $k=0,1,2$ are
(we omit $\int_X$)
\begin{eqnarray}\label{1.32}
V_0^{mfd}(X)&=&c_n,\\
V_2^{mfd}(X)&=& \frac{n}{12}c_n+\frac{1}{6}c_1c_{n-1},\label{1.33},\\
V_4^{mfd}(X)&=& 
\left(\frac{n^2}{48}-\frac{n}{120}\right)\cdot c_n + 
\left((\frac{n}{12}-\frac{1}{30})c_1\right)\cdot c_{n-1}
\label{1.34}\\
&+& \left(\frac{c_2}{10} + \frac{c_1^2}{30}\right)\cdot c_{n-2} 
 + \left(\frac{c_1c_2}{10} - \frac{c_3}{10} - \frac{c_1^3}{30}\right)
\cdot c_{n-3}. \nonumber
\end{eqnarray}

In the case of the projective spaces $\P^n$ the conjectures \ref{t1.2}
are not true for small $k$, see chapter \ref{c7}. It would be interesting
to understand the significance of the Bernoulli moments 
for compact complex manifolds.

\bigskip
When some years ago one of us showed Duco van Straten 
$\Gamma_4^{Ber}(V^{sing},\frac{\alpha_\mu-\alpha_1}{2})$ and the
observation that it is positive in many examples, he conjectured
immediately that there should be a series with signs $(-1)^k$.
We thank him for this idea.

\section{Deformations of higher moments}\label{c2}
\setcounter{equation}{0}

\noindent
Let 
\begin{eqnarray}\label{2.1}
V = \sum_{k=0}^\infty V_{2k}\frac{1}{(2k)!}t^{2k}
\end{eqnarray}
be a formal power series in $t^2$ with variables $V_0,V_2,V_4,...$,
and let $\nu$ be another variable.
We are interested in formal power series
\begin{eqnarray}
\Gamma(V,\nu) &=& \sum_{k=0}^\infty \Gamma_{2k}(V,\nu)\frac{1}{(2k)!}t^{2k}
\quad \hbox{ with }\nonumber\\
\Gamma_{2k}(V,\nu)&\in& \sum_{l=0}^\infty \C[\nu]\cdot V_{2l}\label{2.2}
\end{eqnarray}
which satisfy the following property:
\begin{eqnarray}\label{2.3}
\Gamma(V,\nu)\cdot \Gamma(V',\nu')= \Gamma(V\cdot V',\nu+\nu');
\end{eqnarray}
here $V'$ is a second power series in independent variables,
and $\nu'$ is another variable.

\begin{lemma}\label{t2.1}
A power series $\Gamma(V,\nu)$ as in \eqref{2.2}
satisfies \eqref{2.3} if and only if it takes the form
\begin{eqnarray}\label{2.4}
\Gamma(V,\nu) = \left[\sum_{k=0}^\infty V_{2k}\frac{1}{(2k)!}
\left(\Psi(t)\right)^{k}\right]\cdot \exp(\nu \cdot\Theta (t))
\end{eqnarray}
where
\begin{eqnarray}\label{2.5}
\Psi (t) &=& \sum_{k=1}^\infty \Psi_{2k}\frac{1}{(2k)!}t^{2k},\\
\Theta (t) &=& \sum_{k=1}^\infty \Theta_{2k}\frac{1}{(2k)!}t^{2k},
\label{2.6}
\end{eqnarray}
$\Psi_{2k},\Theta_{2k}\in \C$, or if $\Gamma(V,\nu)=0$.
\end{lemma}

\noindent
{\bf Proof:}
One sees immediately that a power series $\Gamma(V,\nu)$ as in \eqref{2.4}
satisfies \eqref{2.3}. The inverse will be carried out in two steps.

{\bf (I)} We suppose $\Gamma(V,0)=V$ and want to prove
$\Gamma(V,\nu)=V\cdot \exp(\nu\cdot \Theta(t))$. Define
\begin{eqnarray}\label{2.7}
\Phi(t,\nu):= \Gamma(1,\nu)(t)\in \C[\nu][[t]].
\end{eqnarray}
Then $\Phi(t,0)=1$,
\begin{eqnarray}\label{2.8}
\Phi(t,\nu)\cdot \Phi(t,\nu')=\Gamma(1,\nu)\cdot \Gamma(1,\nu') 
= \Gamma(1,\nu+\nu') = \Phi(t,\nu+\nu')
\end{eqnarray}
and
\begin{eqnarray}\label{2.9}
(\log \Phi)(t,\nu) + (\log \Phi)(t,\nu') = (\log \Phi)(t,\nu+\nu').
\end{eqnarray}
One sees easily $(\log \Phi)(t,\nu)\in \C[[\nu]][[t]]$. 
Now $(\log \Phi)(t,0)=\log 1 = 0$ and \eqref{2.9} imply
\begin{eqnarray}\label{2.10}
(\log \Phi)(t,\nu) = \nu\cdot \Theta(t)
\end{eqnarray}
for some $\Theta(t)\in \C[[t]]$. Setting $V'=1$ in \eqref{2.3}
we obtain
\begin{eqnarray}\label{2.12}
\Gamma(V,\nu)\cdot \exp(-\nu\cdot \Theta)=V.
\end{eqnarray}

{\bf (II)} We consider the case with $\nu=0$, that is, without $\nu$.

{\bf Claim:} $\Gamma_0(V,0)=V_0$ or $\Gamma_0(V,0)=0$.

\noindent {\bf Proof:} 
Let $\Gamma_0(V,0)=\lambda_0V_0+...+\lambda_{2l}V_{2l}$ for some $l\geq 0$.
First suppose that $l>0$. Then the special values
$V=V'=1\cdot \frac{1}{(2l)!}t^{2l}$ in \eqref{2.3} yield
\begin{eqnarray}\label{2.13}
\lambda^2_{2l} &=& \Gamma_0(V,0)\cdot \Gamma_0(V',0)\\
&=& \Gamma_0(V\cdot V',0)= \Gamma_0(\frac{1}{((2l)!)^2}t^{4l})=0
\nonumber
\end{eqnarray}
because $4l>2l$. Thus $\lambda_{2l}=0$. Inductively this yields
$\Gamma_0(V,0)=\lambda_0V_0$. Now the same calculation for $l=0$
shows $\lambda^2_0=\lambda_0$, thus $\lambda_0\in \{0;1\}$.
This finishes the proof of the claim.

Now \eqref{2.3} for $V'=1$ gives
\begin{eqnarray}\label{2.14}
\Gamma(V,0)\cdot \Gamma(1,0)=\Gamma(V,0).
\end{eqnarray}
In the case $\Gamma_0=0$ this implies $\Gamma(V,0)=0$. 
We restrict ourselves now to the case $\Gamma_0=V_0$. Then \eqref{2.14}
implies $\Gamma(1,0)=1$. Thus
\begin{eqnarray}\label{2.15}
\Gamma_{2k}(V,0)=\sum_{l>0}\lambda_{kl}\cdot V_{2l}\quad \hbox{ for }
k>0
\end{eqnarray}
is a finite linear combination of terms $V_{2l}$ without the term $V_0$.

Using \eqref{2.15}, we can define $\Psi(t)\in \C[[t^2]]$ by
\begin{eqnarray}\label{2.16}
\Gamma(V_0+V_2\frac{1}{2}t^2)= V_0 + V_2\frac{1}{2}\Psi(t).
\end{eqnarray}
Now we fix $l\in \N$ and 
choose a $V$ with values $V_0=1$ and $V_{2k}=0$ for $k>l$.
As in \cite[Lemma 1.2.1]{Hi} we consider the formal decomposition
of the polynomial $V(t)$ of degree $\leq 2l$,
\begin{eqnarray}\label{2.17}
V(t)=1+\sum_{k=1}^l V_{2k}\frac{1}{(2k)!}t^{2k} = 
\prod_{k=1}^l (1+\beta_{2k}t^2).
\end{eqnarray}
Then
\begin{eqnarray}\nonumber
\Gamma(V(t),0) &=& \prod_{k=1}^l \Gamma(1+\beta_{2k}t^2)
= \prod_{k=1}^l \Gamma(1+\beta_{2k}\Psi(t))\\
&=& 1 + \sum_{k=1}^l V_{2k}\frac{1}{(2k)!}(\Psi(t))^{k}.\label{2.18}
\end{eqnarray}
Because $\Psi(t)$ has no constant term and because the 
$\Gamma_{2k}(V,0)$ are finite linear combinations of the $V_{2k'}$
and because of \eqref{2.15},
this shows for general $V$
\begin{eqnarray}\label{2.19}
\Gamma(V(t),0) = \sum_{k=0}^\infty V_{2k}\frac{1}{(2k)!}(\Psi(t))^{k}.
\end{eqnarray}
This completes the proof. \hfill $\qed$

\begin{remarks}\label{t2.2}
a) The lemma \ref{t2.1} is close to Lemma 1.2.1 and Lemma 1.2.2
in \cite{Hi}. Three differences are the parameter $\nu$ here,
that we do not necessarily have $V_0=1$ and $\Gamma_0(V,\nu)=1$ here
and that here $\Gamma_{2k}(V,\nu)$ is a linear combination of 
the $V_{2l}$, not a quasihomogeneous polynomial.

b) \eqref{2.3} together with the condition $\Gamma(V,0)=V$ restricts
the solutions to the case $\Psi(t)=t^2$. We will only be interested
in this case.
\end{remarks}

\begin{corollary}\label{t2.3}
The Bernoulli moments are characterized by the four properties
\eqref{2.2}, \eqref{2.3}, 
\begin{eqnarray}\label{2.20}
\Gamma_{2k}^{Ber}(V,0)=V_{2k},&&\\
\Gamma_{2k}^{Ber}(V^{sing}(A_\mu),\frac{1}{2})&&
\hbox{ is a polynomial in }w=\frac{1}{\mu+1}\label{2.21}\\
&& \hbox{ for }k\geq 1.\nonumber
\end{eqnarray}
\end{corollary}

\noindent
{\bf Proof:}
The first three conditions show 
$\Gamma^{Ber}(V,\nu)= V\cdot \exp(\nu\cdot \Theta(t))$
for some $\Theta(t)\in \C[[t]]$. By induction on $k>0$,
the condition \eqref{2.21} determines $\Theta_{2k}$ uniquely.
The formulas \eqref{5.4} and \eqref{3.9}
show $\Theta_{2k}=\Theta_{2k}^{Ber}= \frac{-1}{2k}B_{2k}$.
\hfill $\qed$

\bigskip

The following lemma implies lemma \ref{t1.1}.

\begin{lemma}\label{t2.4}
Consider $V(t)\in \R[[t^2]]$ and $\Theta(t)\in \R[[t^2]]$ with 
coefficients $V_{2k}$ and $\Theta_{2k}$ as in \eqref{2.1} and
\eqref{2.6} and $V_0>0$, $-\Theta_2>0$ and $(-1)^k\Theta_{2k}\geq 0$ 
for all $k\geq 2$.
Consider $\Gamma(V,\nu)(t)=V\cdot \exp(\nu\cdot \Theta(t))$.
Fix $k_0\in\N\cup \{\infty\}$.
\begin{list}{}{}
\item[a)] If $k_0\in \N$, there exists a number $\nu\in \R$ such that
\begin{eqnarray}\label{2.22}
\forall\ k\in \N \hbox{ with }k \leq k_0\qquad 
(-1)^k\Gamma_{2k}(V,\nu)\geq 0.
\end{eqnarray}
\item[b)] If a number $\nu\in\R$ satisfies \eqref{2.22} for 
$k_0\in \N\cup\{\infty\}$ then
also any number $\nu'\in \R$ with $\nu'>\nu$ satisfies \eqref{2.22}.
\end{list}
\end{lemma}

\noindent
{\bf Proof:}
a) The polynomial $(-1)^k\Gamma_{2k}(V,\nu)\in \R[\nu]$
has degree $k$. Its term of degree $k$ is
\begin{eqnarray}\label{2.23}
(-1)^kV_0\cdot \Theta_2^k \frac{(2k)!}{k!}\cdot \nu^k.
\end{eqnarray}
It is positive if $\nu>0$, and for large $\nu$ it dominates
$(-1)^k\Gamma_{2k}(V,\nu)$.

b) Consider the two power series $\Theta(it)\in \R[[t^2]]$ 
and $\exp((\nu'-\nu)\cdot \Theta(it))\in \R[[t^2]]$ for some 
fixed $\nu'>\nu$. All their coefficients are nonnegative.
The numbers $(-1)^k\Gamma_{2k}(V,\nu')$ are the coefficients
of 
\begin{eqnarray}\label{2.24}
\Gamma(V,\nu')(it)=\Gamma(V,\nu)(it) 
\cdot \exp((\nu'-\nu)\cdot \Theta(it)).
\end{eqnarray}
If the first $k_0$ coefficients of $\Gamma(V,\nu)(it)$ are nonnegative,
then also the first $k_0$ coefficients 
of $\Gamma(V,\nu')(it)$ are nonnegative. \hfill $\qed$

\begin{remarks}\label{t2.5}
a) In the case of hypersurface singularities, the spectral numbers
satisfy a Thom-Sebastiani property
\cite{Va}\cite{SchS}: Let $f(x_0,...,x_n)$ and $g(y_0,...,y_m)$
be two singularities in different variables with spectral numbers
$\alpha_i$ and $\beta_j$. Then the spectrum of $f+ g$ is the tuple
of numbers
\begin{eqnarray}\label{2.25}
Sp(f+ g)=(\alpha_i+\beta_j+1\ |\ i=1,...,\mu(f);j=1,...,\mu(g)).
\end{eqnarray}
This means that the distribution $\Delta(f+g)(s)$ associated to
$f+g$ (cf. \eqref{1.25}) is the convolution of those associated
to $f$ and $g$,
\begin{eqnarray}\label{2.26}
\Delta(f+g) = \Delta(f)\ast \Delta(g).
\end{eqnarray}
$V^{sing}(f)(2\pi i t)$ is the Fourier transform of $\Delta(f)$. Thus
\begin{eqnarray}\label{2.27}
V^{sing}(f+ g) = V^{sing}(f)\cdot V^{sing}(g).
\end{eqnarray}
With $\nu_1(f):=\alpha_{\mu(f)}(f)-\alpha_1(f)$ and $\nu_2(f):=n+1$,
we get for $j=1,2$ 
\begin{eqnarray}\label{2.28}
&&\Gamma^{Ber}(V^{sing}(f+g),\nu_j(f+g))\\
&=& \Gamma^{Ber}(V^{sing}(f),\nu_j(f)) \cdot
\Gamma^{Ber}(V^{sing}(g),\nu_j(g)).\nonumber
\end{eqnarray}

b) Conjecture \ref{t1.2} (S) [respectively (W)] for a singularity $f$ says 
that all coefficients of 
$\Gamma^{Ber}(V^{sing}(f),\alpha_\mu-\alpha_1)(2\pi it)$ 
[respectively $\Gamma^{Ber}(V^{sing}(f),n+1)(2\pi it)$] are nonnegative 
[respectively positive].
Formula \eqref{2.28} shows part e) of theorem \ref{t1.3}. 

c) Consider a compact complex manifold of dimension $n$ with higher 
moments $V^{mfd}_{2k}(X)$ as in \eqref{1.6} and generating function
$V^{mfd}(X)$. By Serre duality (e.g. \cite[p 123]{Hi}\cite[p 102]{GH}) 
\begin{eqnarray}\label{2.29}
h^{pq}=h^{n-p,n-q}\quad \hbox{ and }\quad \chi_p=\chi_{n-p}.
\end{eqnarray}
Therefore
\begin{eqnarray}\label{2.30}
V^{mfd}(X)= \sum_{p=0}^n \chi_p\cdot \cosh(t(p-\frac{n}{2})
=\sum_{p=0}^n \chi_p \cdot e^{t(p-\frac{n}{2})}.
\end{eqnarray}
If $Y$ is a second compact complex manifold,
the spaces $H^{pq}(X)=H^q(X,\Omega^p)$ and those of $Y$ and $X\times Y$
satisfy the K\"unneth formula (e.g. \cite[p 103]{GH})
\begin{eqnarray}\label{2.31}
H^{\ast,\ast}(X\times Y)\cong H^{\ast,\ast}(X)\otimes H^{\ast,\ast}(Y).
\end{eqnarray}
Therefore 
\begin{eqnarray}\label{2.32}
\chi_p(X\times Y)&=& \sum_{a,b:a+b=p} \chi_a(X)\cdot \chi_b(Y)
\quad \hbox{ and }\\
V^{mfd}(X\times Y)&=&V^{mfd}(X)\cdot V^{mfd}(Y).\label{2.33}
\end{eqnarray}
\end{remarks}

\section{Generalized Bernoulli polynomials}\label{c3}
\setcounter{equation}{0}

\noindent
Define
\begin{eqnarray}\label{3.1}
\Theta^{Ber}(t)= \sum_{k=0}^\infty \Theta^{Ber}_{2k}
\frac{1}{(2k)!}t^{2k}
= \log\frac{t/2}{\sinh (t/2)}.
\end{eqnarray}
The polynomials $A_k(x,\nu)\in \Q[x,\nu]$ for $k\in \N$ are defined by
\begin{eqnarray}\label{3.2}
e^{xt}\cdot \exp\left( \nu\cdot \Theta^{Ber}(t)\right)
= \sum_{k=0}^\infty A_k(x,\nu)\frac{1}{k!}t^{k}.
\end{eqnarray}
They coincide with the classical generalized Bernoulli polynomials
$B^{(\nu)}_k(x)$ up to a shift,
\begin{eqnarray}\label{3.3}
A_k(x,\nu)= B^{(\nu)}_k(x+\frac{\nu}{2}).
\end{eqnarray}
The notation $B^{(\nu)}_k(x)$ was established by N\"orlund 
\cite{No1}\cite{No2}. He and Milne-Thomson \cite{MT} 
studied them systematically for fixed $\nu\in \N$. 
In \cite[p 177]{No1} N\"orlund states that
they had been considered for fixed $\nu\in\N$ already by 
A. Cauchy ($\leq 1890$), E. Lucas (1878), B. Imschenetzky (1883), 
J. Sylvester (1883), D. Sintzof (1890), E. Grigoriew (1898) and 
N. Nielsen (1904). The Bernoulli polynomials $B_k(x)=B_k^{(1)}(x)$ themselves
had first been considered by Jacob Bernoulli ($\leq 1713$), then by Euler.
Since the 19th century the literature on them and on the
Bernoulli numbers $B_k=B_k(0)$ is huge. Their basic properties are 
treated for example in 
\cite{AS}\cite{Er}\cite{Jo}\cite{MT}\cite{No1}\cite{No2}.

In \cite{No1}\cite{No2} there are some remarks about the polynomials
$B_k^{(\nu)}(0)$ in $\nu$. But a study of the $B^{(\nu)}_k(x)$ as polynomials
in $x$ and $\nu$ seems to have been started only in the 60ies, in 
\cite{No3}\cite{No4}\cite{We}. Weinmann \cite{We} seems to be the only one
who shares our point of view that the $A_k(x,\nu)$ have advantages
compared with the $B^{(\nu)}_k(x)$: we have $A_k(-x)=(-1)^kA_k(x)$,
compared to $B_k^{(\nu)}(\nu-x)=(-1)^kB_k^{(\nu)}(x)$, 
and $\deg_\nu A_k(x,\nu)=\left[\frac{k}{2}\right]$, compared to 
$\deg_\nu B_k^{(\nu)}(x)=k$ 
(both are polynomials of degree $k$ in $x$).

The following theorem states well-known or elementary properties
of the Bernoulli numbers and the $A_k(x,\nu)$.

\begin{theorem}\label{t3.1}
a) The Bernoulli numbers satisfy
\begin{eqnarray}\label{3.4}
B_0&=&1,\ B_1=-\frac{1}{2}, \ B_{2k+1}=0 \quad \hbox{ if }k\geq 1,\\
B_{2k} &=& (-1)^{k-1}\frac{2(2k)!}{(2\pi)^{2k}}\zeta(2k)\quad  
\hbox{ if }k\geq 1,\label{3.5}\\
0&=&\sum_{j=0}^{k-1} {k\choose j}B_j \quad \hbox{ if }k\geq 2,\label{3.6}
\end{eqnarray}
\begin{eqnarray}
(B_{2k}\ |\ k=1,...,8) = (\frac{1}{6},-\frac{1}{30},\frac{1}{42},
-\frac{1}{30},\frac{5}{66},-\frac{691}{2730},\frac{7}{6},
-\frac{3617}{510}).\label{3.7}
\end{eqnarray}
\eqref{3.5} shows $B_{2k}=(-1)^{k-1}|B_{2k}|$ if $k\geq 1$ and gives
their asymptotic behaviour because 
$\zeta(s)=\sum_{n=1}^\infty \frac{1}{n^s}\to 1$ fast if $s\to +\infty$.
\eqref{3.6} provides an efficient way to calculate them and shows 
$B_k\in\Q$. The usual definition is via the generating function
\begin{eqnarray}\label{3.8}
\frac{t}{e^t-1}= \sum_{k=0}^\infty B_k\frac{1}{k!}t^k.
\end{eqnarray}
The Bernoulli numbers turn also up in $\Theta^{Ber}(t)$,
\begin{eqnarray}\label{3.9}
\Theta^{Ber}(t) = \sum_{k=1}^\infty \frac{-1}{2k}B_{2k} \frac{1}{(2k)!}t^{2k}.
\end{eqnarray}

b) The polynomials $A_k(x,\nu)$ satisfy
\begin{eqnarray}\label{3.10}
A_k(x,0)&=&x^k,\\
A_{2k+1}(0,\nu)&=&0,\label{3.11}\\
A_{2k}(0,\nu)&\in& (-1)^k\Q_{\geq 0}[\nu], 
\quad \deg_\nu A_{2k}(0,\nu)=k,\label{3.12}
\end{eqnarray}
\begin{eqnarray}
A_k(x,\nu)= \sum_{j=0}^{[k/2]} {k\choose 2j}A_{2j}(0,\nu)\cdot x^{k-2j},
\label{3.13}
\end{eqnarray}
\begin{eqnarray}\label{3.14}
&&A_0(x,\nu)=1,\ A_2(0,\nu)=-\frac{1}{12}\nu,\ 
A_4(0,\nu)=\frac{1}{120}\nu+ \frac{1}{48}\nu^2,\hspace*{1cm}\\
&&A_6(0,\nu)=-(\frac{1}{252}\nu+\frac{1}{96}\nu^2+\frac{5}{576}\nu^3),
\label{3.15}
\end{eqnarray}
\begin{eqnarray}
A_k(x_1+x_2,\nu_1+\nu_2)= \sum_{j=0}^{k} {k\choose j}A_{j}(x_1,\nu_1)
\cdot A_{k-j}(x_2,\nu_2),\hspace*{1cm}
\label{3.16}
\end{eqnarray}
\begin{eqnarray}
A_k(-x,\nu)&=& (-1)^k A_k(x,\nu),\label{3.17}\\
\frac{\paa}{\paa x}A_k(x,\nu)&=&k\cdot A_{k-1}(x,\nu),\label{3.18}\\
\frac{\paa}{\paa \nu}A_k(x,\nu)&=&\sum_{j=1}^{[k/2]} {k\choose 2j}
\frac{-1}{2j}B_{2j} A_{k-2j}(x,\nu),\label{3.19}
\end{eqnarray}
\begin{eqnarray}
A_k(x+\frac{1}{2},\nu+1)- A_k(x-\frac{1}{2},\nu+1)= k\cdot A_{k-1}(x,\nu),
\hspace*{1cm}\label{3.20}\\
\nu\cdot A_k(x\pm\frac{1}{2},\nu+1) = (\nu-k)A_k(x,\nu) + 
k(x\pm \frac{\nu}{2})A_{k-1}(x,\nu),\hspace*{0.5cm}\label{3.21}
\end{eqnarray}
\begin{eqnarray}\label{3.22}
A_k(x,k+1) = \prod_{j=0}^{k-1}(x+\frac{k-1}{2}-j).
\end{eqnarray}
\end{theorem}

\noindent
{\bf Proof:}
a) \eqref{3.8} follows from $B_k=A_k(-\frac{1}{2},1)$ and \eqref{3.2}.
The calculation 
\begin{eqnarray}\label{3.23}
t\frac{\paa}{\paa t}\Theta^{Ber}(t) =1-\frac{1}{2}t\frac{\cosh(t)}{\sinh(t)}
= 1-\frac{1}{2}t - \frac{t}{e^t-1},
\end{eqnarray}
\eqref{3.8}, the fact that $\Theta^{Ber}(t)$ is even and  
$\Theta^{Ber}(0)=0$ show \eqref{3.4} and \eqref{3.9}.

\eqref{3.4} and \eqref{3.17} yield $A_k(-\frac{1}{2},1)=A_k(\frac{1}{2},1)$ 
($=0$ if $k$ is odd and $\geq 3$) for $k\neq 1$. 
Now \eqref{3.10} and \eqref{3.16}
for $(x_1,x_2,\nu_1,\nu_2)=(-\frac{1}{2},1,1,0)$ imply \eqref{3.6}.
With \eqref{3.6} one can calculate \eqref{3.7}. Finally, \eqref{3.5}
is well-known and a consequence of \eqref{3.24}.

b) \eqref{3.10}, \eqref{3.11}, \eqref{3.16}, its special case
\eqref{3.13} and \eqref{3.17} are obvious. 
\eqref{3.12} follows from \eqref{3.2}, \eqref{3.9} and 
$-B_{2k}\in (-1)^k\Q_{>0}$ for $k\geq 1$.
\eqref{3.14} and \eqref{3.15} can be calculated with \eqref{3.19}.
For \eqref{3.18} and \eqref{3.19} one differentiates \eqref{3.2}
by $x$ and $\nu$ and uses \eqref{3.9}. A straightforward calculation
yields \eqref{3.20}. For \eqref{3.21} one applies $t\frac{\paa}{\paa t}$
to \eqref{3.2} and uses \eqref{3.20}. By induction one obtains \eqref{3.22}
from \eqref{3.21} for $k=\nu$.
\hfill $\qed$

\bigskip
Especially interesting for us are the behaviour of $A_k(x,\nu)$ for 
fixed $\nu$ and $k\to \infty$ and the relation to Fourier series. 
Part a) of the following theorem is classical,
part b) is a generalization of a) essentially due to Weinmann \cite{We},
part c) is essentially due to N\"orlund \cite{No4}\cite{No3}.
Part c) contains theorem \ref{t1.4}.
We do not use part b) later, but it fits well to part c) and to 
the conjecture \ref{t1.2} (W).

\begin{theorem}\label{t3.2}
a) Let $f_k:\R\to \R$ be the 1-periodic function with $f_k(x)=A_k(x,1)$ 
for $x\in ]-\frac{1}{2},\frac{1}{2}]$. 
For $k\geq1 $ its Fourier series $\sigma(f_k)$
is 
\begin{eqnarray}\label{3.24}
\sigma(f_{2k})&=& (-1)^{k-1}\frac{2(2k)!}{(2\pi)^{2k}}
\sum_{n=1}^\infty (\frac{(-1)^n}{n^{2k}}\cos(2\pi n x),\\
\sigma(f_{2k-1})&=& (-1)^{k}\frac{2(2k-1)!}{(2\pi)^{2k-1}}
\sum_{n=1}^\infty (\frac{(-1)^n}{n^{2k-1}}\sin(2\pi n x)\label{3.25}.
\end{eqnarray}
For $k\geq 2$ the Fourier series $\sigma(f_k)$ converges uniformly to $f_k$;
the Fourier series $\sigma(f_1)$ converges to $f_1$ on $\R-(\frac{1}{2}+\Z)$.

b) For $\nu\in\N_{\geq 1}$ and $k>\nu$ and 
$x\in [-\frac{\nu}{2},\frac{\nu}{2}]$
\begin{eqnarray}\label{3.26}
A_k(x,\nu) = {k-1\choose \nu-1}\sum_{j=0}^{\nu-1} (-1)^{\nu-1-j}
{\nu-1\choose j}\frac{k}{k-j}\\
\hspace*{4cm}\cdot A_j(x,\nu)\cdot f_{k-j}(x+\frac{\nu-1}{2}).\nonumber
\end{eqnarray}
Replacing the functions $f_{k-j}$ by their Fourier series, one obtains a 
series in $\cos(2\pi n x)$ and $\sin(2\pi nx)$ with polynomial
coefficients, which converges uniformly to $A_k(x,\nu)$ 
on $[-\frac{\nu}{2},\frac{\nu}{2}]$

c) On any compact intervall $I\subset \R$ and for any $\nu\in \R-\Z_{\leq 0}$,
the sequence of polynomials in \eqref{1.21}
tends uniformly to $\cos(2\pi x)$ as $k\to \infty$
and the sequence of polynomials 
\begin{eqnarray}\label{3.27}
(-1)^{k-1} A_{2k-1}(x,\nu)\cdot \frac{(2\pi)^{2k-1}\cdot \Gamma(\nu)}
{2\cdot (2k-1)!\cdot (2k-1)^{\nu-1}}
\end{eqnarray}
tends uniformly to $\sin(2\pi x)$ as $k\to \infty$.
\end{theorem}

{\bf Proof:}
a) See for example \cite[p 37]{Er} for a proof using a contour integral
and \cite[\S 82]{Jo} for a proof in which the Fourier coefficients
are calculated inductively.

b) Weinmann \cite[p 77]{We} generalized the proof of a) via a contour integral.
He obtained the formula which one gets if one replaces the functions
$f_{k-j}$ in \eqref{3.26} by their Fourier series.

We offer a different proof. Suppose that $\nu\in \N_{\geq 1}$ and $k\geq \nu$.
Repeated application of \eqref{3.21} for $x-\frac{1}{2}$ yields the formula
\cite[p 148 (87)]{No2}
\begin{eqnarray}\label{3.28}
A_k(x,\nu) = {k-1\choose \nu-1}\sum_{j=0}^{\nu-1} (-1)^{\nu-1-j}
{\nu-1\choose j}\frac{k}{k-j}\\
\hspace*{4cm}\cdot A_j(x,\nu)\cdot A_{k-j}(x+\frac{\nu-1}{2},1).\nonumber
\end{eqnarray}

{\bf Claim:} The formula remains true if one replaces 
$A_{k-j}(x+\frac{\nu-1}{2},1)$ by $A_{k-j}(x+\frac{\nu-1}{2}-l),1)$
for any $l\in \{0,1,...,\nu-1\}$.

\noindent
{\bf Proof:}
For $l\in \{1,...,\nu-1\}$
the difference of the formulas for $l-1$ and $l$ is, after
dividing by ${k-1\choose \nu-1}/k$,
\begin{eqnarray}\label{3.29}
&&\sum_{j=0}^{\nu-1} (-1)^{\nu-1-j}
{\nu-1\choose j}\frac{1}{k-j}\cdot A_j(x,\nu)\cdot\\ 
&&\hspace*{1cm}
\left(A_{k-j}(x+\frac{\nu-1}{2}-l+1,1)-A_{k-j}(x+\frac{\nu-1}{2}-l,1)\right)
\nonumber\\
&=& \sum_{j=0}^{\nu-1} (-1)^{\nu-1-j}
{\nu-1\choose j}\cdot A_j(x,\nu)\cdot 
A_{k-j-1}(x+\frac{\nu}{2}-l,0)\nonumber\\
&=& \sum_{j=0}^{\nu-1} (-1)^{\nu-1-j}
{\nu-1\choose j}\cdot A_j(x,\nu)\cdot 
(x+\frac{\nu}{2}-l)^{k-j-1} \nonumber\\
&=& (x+\frac{\nu}{2}-l)^{k-\nu}
\sum_{j=0}^{\nu-1} 
{\nu-1\choose j}\cdot A_j(x,\nu)\cdot 
(l-x-\frac{\nu}{2})^{\nu-1-j} \nonumber\\
&=& (x+\frac{\nu}{2}-l)^{k-\nu}
\sum_{j=0}^{\nu-1} 
{\nu-1\choose j}\cdot A_j(x,\nu)\cdot 
A_{\nu-1-j}(l-x-\frac{\nu}{2},0) \nonumber\\
&=& (x+\frac{\nu}{2}-l)^{k-\nu} \cdot A_{\nu-1}(l-\frac{\nu}{2},\nu)
=0.\nonumber
\end{eqnarray}
Here we used \eqref{3.20}, \eqref{3.10}, \eqref{3.16} and \eqref{3.22}.
This shows the claim.

Now for any $x\in[-\frac{\nu}{2},\frac{\nu}{2}]$ there exists an
$l\in \{0,1,...,\nu-1\}$ such that 
$x+\frac{\nu-1}{2}-l\in [-\frac{1}{2},\frac{1}{2}]$.
For $k>\nu$ and $j\leq \nu-1$, the 1-periodic function $f_{k-j}$ 
is continuous and equals $A_{k-j}$ on $[-\frac{1}{2},\frac{1}{2}]$.
Therefore we can replace in \eqref{3.28} 
$A_{k-j}(x+\frac{\nu-1}{2},1)$ by $f_{k-j}(x+\frac{\nu-1}{2})$.

\bigskip
c) Let us fix a compact intervall $I\subset\R$ and a number 
$\nu\in \R-\Z_{\leq 0}$.
It is sufficient to prove that a bound $b>0$ exists such that
for all $k\in\N$ and all $x\in I$
\begin{eqnarray}\label{3.30}
|A_k(x,\nu)\cdot \frac{(2\pi)^{k}\cdot \Gamma(\nu)}
{2\cdot k!\cdot k^{\nu-1}} - \cos(2\pi x-\frac{\pi}{2}k)| 
< b\cdot k^{-7/9}.
\end{eqnarray}
N\"orlund stated this result \cite{No3}\cite{No4}, even with 
$k^{-1}$ instead of $k^{-7/9}$, but for a single $x$
(and with a sign mistake). In \cite{No4} he sketched a proof using \cite{Pe}.
As we could not get \cite{Pe}, we give a proof, following N\"orlund,
but replacing \cite{Pe} by \cite{Er}. 

For any $x\in I$ the function
\begin{eqnarray}\label{3.31}
t\mapsto e^{xt}\exp(\nu\cdot \Theta^{Ber}(t)) =
e^{(x+\frac{1}{2}\nu)t}\left(\frac{t}{e^t-1}\right)^\nu
\end{eqnarray}
is holomorphic on $R:=\C-\{z\in\C\ |\ \Re z=0, \Im z\notin ]-2\pi,2\pi[\}$.
Therefore
\begin{eqnarray}\label{3.32}
\frac{1}{k!} A_k(x,\nu) 
= \frac{1}{2\pi i}
\int_{C_0} t^{-1-k}\cdot
e^{(x+\frac{1}{2}\nu)t}\left(\frac{t}{e^t-1}\right)^\nu\ddd t,
\end{eqnarray}
where $C_0$ is a closed path in $R$ going around $0$ once counterclockwise.
We replace $C_0$ by the union $C_1\cup C_2\cup C_3\cup C_4$ of the
following paths:
$C_1$ is the circle around $2\pi i$ of radius $2\pi k^{-8/9}$,
oriented clockwise, which starts and ends at $2\pi i(1+k^{-8/9})$;
$C_2$ is the half-circle around 0 of radius $2\pi(1+k^{-8/9})$,
oriented counterclockwise, which starts at $2\pi i(1+k^{-8/9})$
and ends at $-2\pi i(1+k^{-8/9})$; $C_3$ and $C_4$ are obtained
from $C_1$ and $C_2$ by the map $\C\to\C$, $z\mapsto -z$.

The purpose of $k^{-8/9}$ is that $(1+k^{-8/9})^k \approx \exp(k^{1/9})$
tends to $\infty$ faster than any power of $k$ if $k\to \infty$,
but that $(1+k^{-16/9})^k\approx \exp(k^{-7/9})\approx 1+O(k^{-7/9})$
tends to 1. The second property 
will allow to replace the function $(1+z)^k$ by the 
function $e^{kz}$ on a disc of radius $k^{-8/9}$ around 0.

We denote by $I_j$, $j=1,2,3,4$ the numbers which are obtained
if one replaces in the right hand side of \eqref{3.32} $C_0$ by $C_j$.
In the following estimate of $|I_2+I_4|$ the factor $t^{-k}$ 
yields the second and the third term 
and $\left(\frac{1}{e^t-1}\right)^\nu$ yields the fourth term; 
\begin{eqnarray}\nonumber
|I_2+I_4|&\leq& const.\cdot (2\pi)^{-k}
\cdot (1+k^{-8/9})^{-k}\cdot k^{8\nu/9}\\
&\leq & const.\cdot (2\pi)^{-k}
\cdot \exp(-k^{1/9})\cdot k^{8\nu/9}.\label{3.33}
\end{eqnarray}
$I_3$ will give the complex conjugate value of $I_1$; so we restrict
ourselves to $I_1$. Let $C_5$ be the circle around 0 of radius
$k^{-8/9}$, oriented counterclockwise, which starts and ends at
$k^{-8/9}$. With the coordinate change $t= 2\pi i(1+\tau)$ we obtain
\begin{eqnarray}\label{3.34}
I_1 = -e^{2\pi i(x+\frac{1}{2}\nu)}
\frac{(2\pi i)^{\nu-k}}{2\pi i}\int_{C_5} e^{(x+\frac{1}{2}\nu)2\pi i\tau}
\frac{(1+\tau)^{\nu-k-1}}{(e^{2\pi i\tau}-1)^{\nu}}\ddd \tau.
\end{eqnarray}
We have for $\tau$ in the disc of radius $k^{-8/9}$ around 0
\begin{eqnarray}\label{3.35}
e^{(x+\frac{1}{2}\nu)2\pi i\tau}(1+\tau)^{\nu-1}\approx 1+O(k^{-8/9}),\\
(1+\tau)^{-k}\approx e^{-\tau k}\cdot (1+ O(k^{-7/9})),\label{3.36}\\
(e^{2\pi i\tau}-1)^{-\nu}\approx (2\pi i\tau)^{-\nu}\cdot (1+O(k^{-8/9})).
\label{3.37}
\end{eqnarray}
Therefore
\begin{eqnarray}\label{3.38}
I_1 = -e^{2\pi i(x+\frac{1}{2}\nu)}
\frac{(2\pi i)^{-k}}{2\pi i}
\int_{C_5} \frac{e^{-\tau k}}{\tau^\nu}(1+k^{-7/9}g_k(x,\tau))
\ddd \tau,
\end{eqnarray}
where $g_k(x,\tau):I\times \{z\ |\ |z|\leq k^{-8/9}\}\to \C$ is real
analytic in $x$ and holomorphic in $z$ and bounded independently of 
$k,x,z$. Formula (6) in \cite[p 14]{Er} says
\begin{eqnarray}\label{3.39}
- e^{\pi i\nu}\frac{1}{2\pi i}
\int_{C_6} \frac{e^{-\tau k}}{\tau^\nu}\ddd \tau
=\frac{k^{\nu-1}}{\Gamma(\nu)},
\end{eqnarray}
where $C_6$ is a path from $+\infty$ to $+\infty$ circulating once
counterclockwise around $0$. Therefore
\begin{eqnarray}\label{3.40}
I_1 &=& e^{2\pi ix-\frac{\pi}{2}ik}\frac{k^{\nu-1}}{(2\pi)^k\Gamma(\nu)}\\
&-&k^{-7/9}e^{2\pi i(x+\frac{1}{2}\nu)}
\frac{(2\pi i)^{-k}}{2\pi i}
\int_{C_5} \frac{e^{-\tau k}}{\tau^\nu}g_k(x,\tau)
\ddd \tau,\label{3.41}\\
&-& e^{2\pi i(x+\frac{1}{2}\nu)}
\frac{(2\pi i)^{-k}}{2\pi i}
\int_{C_6-C_5} \frac{e^{-\tau k}}{\tau^\nu}
\ddd \tau \label{3.42}.
\end{eqnarray}
The integral \eqref{3.42} can be estimated easily. Its vanishing order
is dominated by $(2\pi)^{-k}\cdot\exp(-k^{1/9})$.
In order to estimate the integral \eqref{3.41}, we replace $C_5$
by $(-C_7)\cup C_8 \cup C_7$, where $C_7$ is the straight line
from $k^{-1}$ to $k^{-8/9}$ and $C_8$ is the circle around 0
of radius $k^{-1}$, oriented counterclockwise, which starts and
ends at $k^{-1}$. With the coordinate change $\widetilde \tau =k\tau$, 
it is easy to see that for $j=7,8$ the integral
\begin{eqnarray}\label{3.43}
k^{1-\nu}\cdot\int_{C_j} 
\left|\frac{e^{-\tau k}}{\tau^{\nu}}\right|\ddd \tau
\end{eqnarray}
is bounded independently of $k$. Therefore \eqref{3.41}
is of order $k^{-7/9}\cdot k^{\nu-1}/(2\pi)^k$.

For $I_3$ we get the complex conjugate result. Thus
\begin{eqnarray}\label{3.44}
\frac{(2\pi)^k\Gamma(\nu)}{2\cdot k^{\nu-1}}(I_1+I_2+I_3+I_4) = 
\cos(2\pi x-\frac{\pi}{2}k) + O(k^{-7/9}).
\end{eqnarray}
This finishes the proof. \hfill $\qed$

\begin{remarks}\label{t3.3}
a) The asymptotic behaviour of the polynomials $A_k(x,\nu)$ has also
been studied in \cite{We} and \cite{No4} in the case $k=k_0+r$ 
and $\nu=\nu_0+r$ with $r\to \infty$. N\"orlund obtains that
a suitable normalization of $B_k^{(\nu)}(x)$ tends to 
$\frac{1}{\Gamma(1-x)}$, Weinmann finds that a suitable normalization
of $A_k(x,\nu)$ tends to a linear combination of $\cos(\pi x)$ and 
$\sin(\pi x)$ with polynomial coefficients in $r$.
So both find the average intervall 1 between neighbouring zeros
of $A_k(x,\nu)$. In the case $A_k(x,k+1)$ this is obvious because of 
formula \eqref{3.22}. In theorem \ref{t3.2} c) we had $\frac{1}{2}$.

b) If $k$ and $\nu$ tend to $\infty$ with a larger [or smaller] 
fixed quotient $\nu/k$, we expect a larger [or smaller] 
average intervall between neighbouring zeros of $A_k(x,\nu)$. 

c) Also, we expect that $A_k(x,\nu)$ has the 
maximal number $k$ of zeros if $\nu\geq k$. This is clear
for $\nu=k+1$ by \eqref{3.22}. Because of \eqref{3.18}
it holds for all $\nu\in \N$ with $\nu\geq k+1$: for these $\nu$ 
\begin{eqnarray}\label{3.45}
A_k(x,\nu)= \frac{(\nu-1)!}{k!}
\frac{\paa^{\nu-1-k}}{\paa x^{\nu-1-k}}A_{\nu-1}(x,\nu).
\end{eqnarray}

d) The (real) zeros of the Bernoulli polynomials and thus of the 
polynomials $A_k(x,1)$ are well understood 
(\cite{De1}\cite{De2}\cite{In}\cite{Le} and references there).
Inkeri \cite{In} showed that the number of zeros of the
Bernoulli polynomials and of the polynomials $A_k(x,1)$ tends
to $\frac{2k}{\pi e}$ as $k\to \infty$. His results are much more precise. 
Delange \cite{De1}\cite{De2} even refined Inkeri's results to such 
a precision that he can derive without effort that
$A_{1000000}(x,1)$ has $234204$ zeros.
Also the positions of the zeros are well understood.

e) If $k$ is fixed and $\nu$ tends to $\infty$, then the zeros of
$A_k(x,\nu)$ tend to $\sqrt{\nu}\cdot c_j$, $j=1,...,k$
with $c_1\leq ...\leq c_k$. This follows from \eqref{3.12}
and \eqref{3.13}. We expect that the numbers $c_1,...,c_k$ are
all different.
So for large $\nu$ the polynomial $A_k(x,\nu)$ is oscillating
around $0$ only for $|x|\leq c_k\cdot \sqrt{\nu}$. For the conjectures
\ref{t1.2} the intervall $[-\frac{\nu}{2},\frac{\nu}{2}]$ is relevant.
\end{remarks}

We conclude with a discussion of $A_2(x,\nu)$ and $A_4(x,\nu)$.

\begin{examples}\label{t3.4}
a) The polynomial $-A_2(x,\nu)=-x^2+\frac{1}{12}\nu$ has the zeros
$\pm \sqrt{\frac{1}{12}\nu}$. The positive zero is smaller than 
$\frac{\nu}{2}$ if $\nu>\frac{1}{3}$.

b) The polynomial $A_4(x,\nu)= x^4-\frac{\nu}{2} x^2
+(\frac{\nu}{120}+\frac{\nu^2}{48})$ 
has two minima at $\pm x_0=\pm\sqrt{\frac{\nu}{4}}$ and a local
maximum at 0. It has four zeros if $\nu>\frac{1}{10}$.
If $\nu>1$ then $x_0<\frac{\nu}{2}$. 
If $\nu>1,768$ then the largest zero is smaller than 
$\frac{\nu}{2}$.
For large $\nu$ the positive zeros are approximately 
$x_0(1\pm \sqrt{\frac{2}{3}})=x_0(1\pm 0,8165)$.
\end{examples}

\section{Interpretation}\label{c4}
\setcounter{equation}{0}

\noindent
The conjectures \ref{t1.2} are about the higher moments of
the spectral numbers of a singularity. Nevertheless it is 
difficult to derive from them concrete information on the
distribution of the spectral numbers. 
The following remarks point to different aspects of this problem.

\begin{remarks}\label{t4.1}
a) The meaning of the conjecture \eqref{1.3}, that is, the case $k=1$,
is clear: the variance is bounded from above. Also for $k\to \infty$
the meaning of the conjectures \ref{t1.2} is clear: by the discussion
after theorem \ref{t1.4} they boil down to the topological 
statement that the sign of the trace of the monodromy is $(-1)^{n-1}$.
But for $k=2$ and any fixed $k\geq 2$ the meaning of 
the conjectures \ref{t1.2} is not at all clear.
If $k$ is small compared to $\nu$, then by remark \ref{t3.3} e)
the polynomial $(-1)^kA_{2k}(x,\nu)$ is oscillating around 0 only 
for $const.\cdot\sqrt{\nu}$ and has the sign $(-1)^k$ outside, 
whereas the conjectures \ref{t1.2} are concerned with the whole intervall
$[-\frac{\nu}{2},\frac{\nu}{2}]$.

b) Because of \eqref{1.21} and \eqref{1.20}, the power series 
$\Gamma^{Ber}(V^{sing}(f),\nu)(2\pi it)$ has 
the radius of convergence 1 for any $\nu\in \R-\Z_{\leq 0}$. 
The conjecture \ref{t1.2} (W) [respectively (S)]
says that all coefficients are positive [nonnegative]
if $\nu=n+1$ [$\nu=\alpha_\mu-\alpha_1$]. What does this say about the
function?

c) It would be good to establish an inverse Fourier transform 
$\FF^{(-1)}(f)(s)$ of the function $\Gamma^{Ber}(V^{sing}(f),n+1)(2\pi it)$.
Then \eqref{1.29} could be rewritten as
\begin{eqnarray}\label{4.1}
\Delta(f)(s) = (\FF^{(-1)}(f)\ast \Delta^{(n+1)})(s);
\end{eqnarray}
this could help to give a better answer to K. Saito's hope 
\cite[p 202, (2.5) ii)]{SK2} that the limit distribution $\Delta^{(n+1)}(s)$
should be a bound of the distributions $\Delta(f)(s)$ of the spectral
numbers of singularities $f$.

d) K. Saito formulated some questions connected with this hope 
\cite[p 203, (2.8)]{SK2}: Is 
\begin{eqnarray}\label{4.2}
|\{j\ |\ \alpha_j\leq -\frac{1}{2}\}| &<& \frac{\mu}{(n+1)!2^{n+1}}\ ?\\
|\{j\ |\ \alpha_j< 0\}| &<& \frac{\mu}{(n+1)!}\ ?\label{4.3}
\end{eqnarray}
For $n=2$ a yes to the second question (with $\alpha_j\leq 0$ instead of
$\alpha_j<0$) is equivalent \cite{SM1} to Durfee's 
conjecture \cite{Du} that the geometric genus of a singularity is $<\mu/6$.

But the conjectures \ref{t1.2} do not answer these questions, see the 
example \ref{t4.2}. They give only weaker inequalities. 
If one could combine them with (unknown) 
statements about ``series'' of spectral numbers, 
they might give stronger estimates.

e) The conjectures \ref{t1.2} point to relations which should be explored
and structures which have yet to be established. On the one hand
there is the similarity of $V^{sing}(f)$ and $V^{mfd}(X)$ 
for compact complex manifolds $X$. Could one hope to
establish for singularities some of the central characters
in chapter \ref{c7}, Chern classes and Hirzebruch-Riemann-Roch theorem?

On the other hand, the conjecture \eqref{1.3} was found \cite{He1}\cite{He2}
by looking at the G-function of Frobenius manifolds \cite{DZ}\cite{Gi}.
In the singularity case, this is a distinguished 
holomorphic function on the Frobenius manifold, that is, the base space
of a semiuniversal unfolding.
Its derivative by the Euler field is just the constant
$-\frac{1}{4}\cdot \Gamma^{Ber}_2(V^{sing}(f),\alpha_\mu-\alpha_1)$.
In the quantum cohomology case, the G-function is the generating
function of the genus 1 Gromov-Witten invariants (the generating
function of the genus 0 invariants gives the Frobenius manifold).
In that case one has generating functions for the invariants of all
genera. Are they related to the higher Bernoulli moments?

These two structures, Chern classes and Frobenius manifolds, 
might also have a chance to provide techniques for proving the 
conjectures \ref{t1.2} in general.
\end{remarks}

\begin{example}\label{t4.2}
The conjecture \ref{t1.2} (W) does not imply the inequality 
\eqref{4.3} in the case $n=2$. We consider an abstract spectrum
with spectral numbers $-\frac{1}{2}$, $0$, and $\frac{1}{2}$,
with multiplicities $r$, $\mu-2r$, $r$, where $0\leq r \leq \frac{\mu}{2}$,
$r\in \R$. Then 
\begin{eqnarray}\label{4.4}
\Gamma^{Ber}_{2k}(V,3)=2r A_{2k}(\frac{1}{2},3)+(\mu-2r)A_{2k}(0,3).
\end{eqnarray}
For $k=1$ $A_{2}(\frac{1}{2},3)=0$; so it does not give any
restriction on $r$. For large $k$ 
$A_{2k}(\frac{1}{2},3)\approx -A_{2k}(0,3)$ by theorem \ref{t1.4};
this gives in the limit the restriction $2r\leq \mu-2r$,
that is, $r\leq \frac{\mu}{4}$, and not $r\leq \frac{\mu}{6}$.
\end{example}

\section{Quasihomogeneous singularities}\label{c5}
\setcounter{equation}{0}

\noindent
A quasihomogeneous singularity $f(x_0,...,x_n)$ has unique (up to 
ordering) normalized weights $w_0,...,w_n\in \Q\cap ]0,\frac{1}{2}]$ such that 
$f$ has weighted degree 1 \cite{SK1}. We will always use these weights.

The starting point of the formulas in this chapter is the following 
well known generating function of the spectrum $\alpha_1,...,\alpha_\mu$ 
of a quasihomogeneous singularity:
\begin{eqnarray}\label{5.1}
\sum_{j=1}^\mu T^{\alpha_j-\frac{n-1}{2}} 
=\prod_{i=0}^n \frac{T^{w_i-\frac{1}{2}}-T^{\frac{1}{2}}}{1-T^{w_i}}.
\end{eqnarray}
Because of \eqref{1.18}, $V^{sing}(f)$ is given by the following 
formula, interpreted as a formal power series in $t$.
\begin{eqnarray}\label{5.2}
V^{sing}(f)= \prod_{i=0}^n \frac{e^{(w_i-\frac{1}{2})t}-e^{\frac{1}{2}t}}
{1-e^{w_it}}.
\end{eqnarray}
The proofs of theorem \ref{t5.1} to theorem \ref{t5.4} will be given
after theorem \ref{t5.4}. The Bernoulli numbers $B_{2k}$ satisfy
$B_{2k}\in(-1)^{k-1}\Q_{>0}$ for $k\geq 1$ and $B_0=1$ (theorem \ref{t3.1}).

\begin{theorem}\label{t5.1}
Let $f(x_0,...,x_n)$ be a quasihomogeneous singularity with 
normalized weights $w_0,...,w_n$. Then
\begin{eqnarray}\label{5.3}
V^{sing}(f) = \prod_{i=0}^{n}\left[ \sum_{k=0}^\infty 
 \left( w_i^{2k}\frac{2}{2k+1}B_{2k+1}(\frac{1}{2w_i})\right)
\frac{1}{(2k)!}t^{2k} \right],\\
\Gamma^{Ber}(V^{sing}(f),n+1)
= \prod_{i=0}^{n} \left[ \sum_{k=0}^\infty (-B_{2k})(1-w_i^{2k-1})
\frac{1}{(2k)!}t^{2k}\right]. 
\label{5.4}
\end{eqnarray}
\eqref{5.4} shows conjecture \ref{t1.2} (W) for $f$, as 
$(-B_{2k})(1-w_i^{2k-1})$ has the sign $(-1)^k$ for any $k\geq 0$.
\end{theorem}

The calculation \eqref{5.22} of the formula \eqref{5.4}
will also be useful for the conjecture \ref{t1.2} (W)
in the case of curve singularities.

\begin{theorem}\label{t5.2}
Conjecture \ref{t1.2} (S) is true for the hyperbolic 
singularities $T_{pqr}$. Then $\alpha_\mu-\alpha_1=1$ and 
\begin{eqnarray}\label{5.5}
\Gamma^{Ber}(V^{sing}(T_{pqr}),1) &= &
\sum_{k=0}^\infty \frac{1}{(2k)!}t^{2k}\cdot \\
&&\left[ B_{2k}\cdot 
\left(-1+\frac{1}{p^{2k-1}}+\frac{1}{q^{2k-1}}+\frac{1}{r^{2k-1}}\right)
\right]. \nonumber
\end{eqnarray}
\end{theorem}

\begin{proposition}\label{t5.3}
Define $Q(t,w)\in \Q[w][[t^2]]$ by 
\begin{eqnarray}\label{5.6}
Q(t,w) = \frac{w}{1-w}  \left(
\frac{e^{(w-\frac{1}{2})t}-e^{\frac{1}{2}t}}{1-e^{wt}}\right)
\exp((1-2w)\Theta^{Ber}(t)).
\end{eqnarray}
a) Then
\begin{eqnarray}
&&Q(t,w)\nonumber\\
&=&\exp\left(\Theta^{Ber}(wt)-\Theta^{Ber}((1-w)t)
+(1-2w)\Theta^{Ber}(t)\right)\label{5.7}\\
&=& \exp\left(\sum_{k=1}^\infty \frac{-1}{2k}B_{2k}p_{2k}(w)
\frac{1}{(2k)!}t^{2k}\right),\label{5.8}
\end{eqnarray}
where
\begin{eqnarray}\label{5.9}
p_{2k}(w) &=& 1-2w + w^{2k}-(1-w)^{2k}.
\end{eqnarray}
b) The first three of the polynomials $p_{2k}$ are 
\begin{eqnarray}\label{5.10}
p_2(w)&=&0,\\
p_4(w)&=&4(\frac{1}{2}-w)w(1-w),\label{5.11}\\
p_6(w)&=&6(\frac{1}{2}-w)w(1-w)(\frac{4}{3}-(w(1-w)).\label{5.12}
\end{eqnarray}
For $k\geq 2$, the polynomial $p_{2k}$ has three simple zeros at 
$0,\frac{1}{2},1$ and no other zeros. It is negative for 
$w\in ]-\infty,0[\cup ]\frac{1}{2},1[$ and positive for 
$w\in ]0,\frac{1}{2}[\cup ]1,+\infty[$. 

c) The polynomials $Q_{2k}(w)$ in 
$Q(t,w)=\sum_kQ_{2k}(w)\frac{1}{(2k)!}t^{2k}$ satisfy
\begin{eqnarray}\label{5.13}
Q_0=1,\quad Q_2=0,\quad Q_4=\frac{1}{30}\cdot\frac{1}{4}p_4,
\quad Q_6=-\frac{1}{42}\cdot\frac{1}{6}p_6,
\end{eqnarray}
and for $k\geq 2$
\begin{eqnarray}\label{5.14}
(-1)^kQ_{2k}(w)>0 \quad \hbox{ if }w\in ]0,\frac{1}{2}[\cup ]1,+\infty[.
\end{eqnarray}
They have simple zeros at $0,\frac{1}{2},1$.
\end{proposition}

We expect that they also satisfy 
\begin{eqnarray}\label{5.15}
&& (-1)^kQ_{2k}(w)<0 \quad \hbox{ if }w\in ]-\infty,0[\cup ]\frac{1}{2},1[,
\end{eqnarray}
but we do not have a proof.

\begin{theorem}\label{t5.4}
Let $f(x_0,...,x_n)$ be a quasihomogeneous singularity with 
normalized weights $w_0,...,w_n$. Then
\begin{eqnarray}\label{5.16}
\Gamma^{Ber}(V^{sing}(f),\alpha_\mu-\alpha_1)
= \mu\prod_{i=0}^n Q(t,w_i).
\end{eqnarray}
\eqref{5.16} and \eqref{5.14} show conjecture \ref{t1.2} (S) for $f$.
\eqref{5.10} -- \eqref{5.13} show
\begin{eqnarray}\label{5.17}
\Gamma_2^{Ber}(V^{sing}(f),\alpha_\mu-\alpha_1)
&=& 0,\\
\Gamma_4^{Ber}(V^{sing}(f),\alpha_\mu-\alpha_1)
&=& \frac{1}{30} \mu \sum_{i=0}^{n} (\frac{1}{2}-w_i)w_i(1-w_i),
\label{5.18}\\
\Gamma_6^{Ber}(V^{sing}(f),\alpha_\mu-\alpha_1)
&=& \frac{1}{42} \mu \sum_{i=0}^{n} (\frac{1}{2}-w_i)w_i(1-w_i)\cdot\nonumber\\
&&\cdot(w_i(1-w_i)-\frac{4}{3}). \label{5.19}
\end{eqnarray}
\end{theorem}

\eqref{5.17} says that in the case of a quasihomogeneous
singularity one has equality in \eqref{1.3}. 
The first proof in \cite{He1}\cite{He2} used Frobenius manifolds,
the second proof in \cite{Di} was elementary and used the formula \eqref{5.1}.
The third proof here in chapter \ref{c5} also uses this formula.
But it is more general and yields also the other formulas 
in the theorems \ref{t5.1} to \ref{t5.4}.

$Q_2=0$ is responsable for \eqref{5.17} and for the simplicity of the formulas 
\eqref{5.18} and \eqref{5.19}.
For $k\geq 4$ one has also products of the $Q_{2l}(t,w_i)$ in the formulas
for the Bernoulli moments $\Gamma_{2k}^{Ber}(V^{sing}(f),\alpha_\mu-\alpha_1)$.

\bigskip\noindent
{\bf Proof of theorem \ref{t5.1}:}
One derives from \eqref{3.1} -- \eqref{3.3} the classical generating
function for the Bernoulli polynomials 
$B_k(x)=B_k^{(1)}(x)=A_k(x-\frac{1}{2},1)$:
\begin{eqnarray}\label{5.20}
\frac{te^{xt}}{e^t-1}= \sum_{k=0}^\infty B_k(x)\frac{1}{k!}t^k.
\end{eqnarray}
The following calculation shows \eqref{5.3}.
\begin{eqnarray}\label{5.21}
&&\sum_{k=0}^\infty
\left( w^{2k}\frac{2}{2k+1}B_{2k+1}(\frac{1}{2w})\right)
\frac{1}{(2k)!}t^{2k}\\
&=& \frac{2}{wt}\sum_{k=0}^\infty B_{2k+1}(\frac{1}{2w})
\frac{1}{(2k+1)!}(wt)^{2k+1}\nonumber\\
&=& \frac{1}{wt}\left( \frac{wte^{\frac{1}{2w}wt}}{e^{wt}-1}
-\frac{-wte^{\frac{1}{2w}(-wt)}}{e^{-wt}-1}\right)\nonumber\\
&=& \frac{e^{\frac{1}{2}t}}{e^{wt}-1} 
+ \frac{e^{-\frac{1}{2}t}}{e^{-tw}-1}
= \frac{-e^{\frac{1}{2}t} + e^{(w-\frac{1}{2})t}}{1-e^{wt}}.\nonumber
\end{eqnarray}
The coefficient of $\frac{1}{(2k)!}t^{2k}$ in the first line of 
\eqref{5.21} is not a polynomial in $w$, but has a pole of order
1 at $w=0$. The calculation \eqref{5.22} shows that the multiplication
by $\exp(\frac{1}{2}\Theta^{Ber}(t))$ cancels these poles for $k\geq 1$.
The coefficients $\Theta^{Ber}_{2k}=-\frac{1}{2k}B_{2k}$ 
are inductively determined by this property. 
This explains the characterisation of the Bernoulli moments 
in corollary \ref{t2.3}.

Formula \eqref{5.4} is a consequence of \eqref{5.2} and the 
following calculation, which uses at the end \eqref{3.8} and 
$B_1=-\frac{1}{2}$.
\begin{eqnarray}\label{5.22}
&& \frac{e^{(w-\frac{1}{2})t}-e^{\frac{1}{2}t}}{1-e^{wt}}
\cdot \exp(\Theta^{Ber}(t))\\
&=& \frac{e^{(w-\frac{1}{2})t}-e^{\frac{1}{2}t}}{1-e^{wt}}
\cdot \frac{t\cdot e^{\frac{1}{2}t}}{e^t-1} 
= \frac{t e^{wt}-te^t}{(1-e^{wt})(e^t-1)} \nonumber
\end{eqnarray}
\begin{eqnarray}
&=& -\frac{t}{e^t-1} + \frac{t}{e^{wt}-1}\nonumber\\
&=& -\left(\frac{t}{e^t-1} +\frac{1}{2}t\right) + 
\frac{1}{w}\left(\frac{wt}{e^{wt}-1}+\frac{1}{2}wt\right) \nonumber\\
&=& \sum_{k=0}^\infty (-B_{2k})(1-w^{2k-1})\frac{1}{(2k)!}t^{2k}.\nonumber
\end{eqnarray}
\hfill $\qed$

\bigskip
\noindent
{\bf Proof of theorem \ref{t5.2}:}
The generating function of the spectrum of the hyperbolic 
surface singularity $T_{pqr}$ is
\begin{eqnarray}\label{5.23}
\sum_{j=1}^\mu T^{\alpha_j} =T^0 + T^1 + 
\frac{T^{1/p}-T}{1-T^{1/p}} + \frac{T^{1/q}-T}{1-T^{1/q}} + 
\frac{T^{1/r}-T}{1-T^{1/r}}.
\end{eqnarray}
Because of \eqref{1.18}
\begin{eqnarray}\label{5.24}
V^{sing}(T_{pqr})= e^{-\frac{1}{2}t}\left(1+e^t+
\frac{e^{\frac{1}{p}t}-e^t}{1-e^{\frac{1}{p}t}} + 
\frac{e^{\frac{1}{q}t}-e^t}{1-e^{\frac{1}{q}t}} + 
\frac{e^{\frac{1}{r}t}-e^t}{1-e^{\frac{1}{r}t}}\right).
\end{eqnarray}
Then, using \eqref{5.22} for $w=\frac{1}{p},\frac{1}{q},\frac{1}{r}$, 
one finds
\begin{eqnarray}\label{5.25}
&&\Gamma^{Ber}(V^{sing}(T_{pqr}),1) \\
&=& \left(e^{-\frac{1}{2}t}+e^{\frac{1}{2}t}\right)\
\frac{te^{\frac{1}{2}t}}{e^t-1} + 
\left(-\frac{t}{e^t-1} + \frac{t}{e^{\frac{1}{p}t}-1}\right) + ... \nonumber\\
&=& \left(2\frac{t}{e^t-1}+t\right) + 
\left(-\frac{t}{e^t-1} + \frac{t}{e^{\frac{1}{p}t}-1}\right) + ... \nonumber\\
&=& \sum_{k=0}^\infty B_{2k}
\left(-1+\frac{1}{p^{2k-1}}+\frac{1}{q^{2k-1}}+\frac{1}{r^{2k-1}}\right)
\frac{1}{(2k)!}t^{2k} .\nonumber
\end{eqnarray}
\hfill $\qed$

\bigskip
\noindent
{\bf Proof of proposition \ref{t5.3}:}
a) \eqref{5.7} follows from
\begin{eqnarray}\label{5.26}
\frac{w}{1-w}  \left(
\frac{e^{(w-\frac{1}{2})t}-e^{\frac{1}{2}t}}{1-e^{wt}}\right)
= \frac{\frac{1}{2}wt}{\frac{1}{2}(1-w)t}
\frac{\sinh(\frac{1}{2}(1-w)t)}{\sinh(\frac{1}{2}wt)}.
\end{eqnarray}
and the definition \eqref{3.1} of $\Theta(t)$.
\eqref{5.8} follows from \eqref{3.9}.

b) For $k\geq 2$ one calculates 
\begin{eqnarray}\label{5.27}
&& p_{2k}(0)=p_{2k}(\frac{1}{2})=p_{2k}(1)=0,\\
&& p_{2k}'(0)=p_{2k}'(1)=2k-2>0,\label{5.28} \\
&& p_{2k}'(\frac{1}{2})=-2+k\cdot2^{3-2k}<0,\label{5.29}\\
&& p_{2k}'''(w)=2k(2k-1)(2k-2)(w^{2k-3}+(1-w)^{2k-3})>0\hspace*{0.5cm}
\label{5.30}.
\end{eqnarray}
Because of \eqref{5.30} the simple zeros of $p_{2k}$ at $0,\frac{1}{2},1$
are the only zeros of $p_{2k}$ for $k\geq 2$.

c) \eqref{5.13} and \eqref{5.14} follow from a) and b). 
The $Q_{2k}$ have simple zeros at $0,\frac{1}{2},1$, because a
calculation shows for $k\geq 2$
\begin{eqnarray}\label{5.31}
Q_{2k}'(0)= Q_{2k}'(1) &=& -B_{2k}(1-\frac{1}{k}),\\
Q_{2k}'(\frac{1}{2})&=& B_{2k}(\frac{1}{k}-\frac{1}{2^{2k-2}}).\label{5.32}
\end{eqnarray}
\hfill $\qed$

\bigskip
\noindent
{\bf Proof of theorem \ref{t5.4}:}
\eqref{5.16} follows from \eqref{5.6}, \eqref{5.2} and 
\begin{eqnarray}\label{5.33}
\alpha_\mu-\alpha_1=\sum_{i=0}^n (1-2w_i),
\quad \mu=\prod_{i=0}^n \left(\frac{1}{w_i}-1\right).
\end{eqnarray}
The rest is a consequence of proposition \ref{t5.3}.
\hfill $\qed$

\section{Curve singularities}\label{c6}
\setcounter{equation}{0}

\begin{theorem}\label{t6.1}
Conjecture \ref{t1.2} (W) is true for any irreducible
curve singularity.
\end{theorem}

\begin{figure}[H]
\centerline{\input{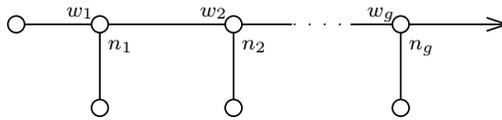}}
\caption{Eisenbud and Neumann diagram of an irreducible germ of a curve}
\label{Figure : EN diagram of an irreducible germ}
\end{figure}

\noindent
{\bf Proof:} 
Suppose that the Puiseux pairs of the irreducible germ of curve $f$ are 
$(n_1,r_1), \ldots, (n_g,r_g)$. Then with $w_1=r_1$, and for 
$k \geq 1$, $w_{k+1}=r_{k+1}-r_kn_{k+1}+n_kn_{k+1}w_k$, the Eisenbud and 
Neumann diagram is given by figure 
\ref{Figure : EN diagram of an irreducible germ} 
(see \cite{Ne} for a rapid overview).
Furthermore let us introduce $n'_k=n_{k+1} \ldots n_g$ for $1 \leq k 
\leq g-1$ and $n'_g=1$.

\begin{figure}[H]
\centerline{\input{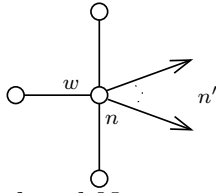}}
\caption{Eisenbud and Neumann diagram of a quasihomogeneous isolated 
curve singularity}
\label{Figure : EN diagram of an quasihomogeneous singularity of curve}
\end{figure}

From \cite{Br3}, we have a formal decomposition of this diagram 
in terms of the Newton non-degenerate and commode germs.
If we denote by $D(w,n,n')$ the diagram given by figure \ref{Figure : 
EN diagram of an quasihomogeneous singularity of curve}, 
where $w$, $n$ are coprime positive integers and $n'$ is a positive integer 
then the decomposition is
\begin{eqnarray}\label{6.1}
D(w_1,n_1,n'_1)+ \sum_{k=1}^{g-1} 
\left( D(w_{k+1},n_{k+1},n'_{k+1})-D(w_kn_k,1,n'_k) \right).
\end{eqnarray}
This gives
\begin{multline}\label{6.2}
Sp(f)= Sp(D(w_1,n_1,n'_1)) \\
+\sum_{k=1}^{g-1} \left( Sp(D(w_{k+1},n_{k+1},n'_{k+1}))-Sp(D(w_kn_k,1,n'_k)) 
\right) .
\end{multline}
More precisely, the generating function $\sum_{i=1}^\mu T^{\alpha_i+1}$ is 
\begin{multline}\label{6.3}
\frac{T^{\frac{1}{n'_0}}-T}{1-T^{\frac{1}{n'_0}}} \cdot 
        \frac{T^{\frac{1}{w_1n'_1}}-T}{1-T^{\frac{1}{w_1n'_1}}}\\
+\sum_{k=1}^{g-1}
\left(
 \frac{T^{\frac{1}{w_{k+1}n'_{k+1}}}-T}    {1-T^{\frac{1}{w_{k+1}n'_{k+1}}}}
- \frac{T^{\frac{1}{w_kn'_{k-1}}}-T}{1-T^{\frac{1}{w_kn'_{k-1}}}} 
\right)
\frac{T^{\frac{1}{n'_k}}-T}{1-T^{\frac{1}{n'_k}}}. 
\end{multline}
From the quasihomogeneous case (the calculation \eqref{5.22} and the formula
\eqref{5.4}), we know that the first term verifies
the conjecture (W). Now to prove the conjecture (W), 
it is sufficient to prove it for 
\begin{eqnarray}\label{6.4}
\left(
\frac{T^{\frac{1}{w_2}}-T}{1-T^{\frac{1}{w_2}}}
-
\frac{T^{\frac{1}{w_1n_1n_2}}-T}{1-T^{\frac{1}{w_1n_1n_2}}}  
\right)
\frac{T^{\frac{1}{n_2}}-T}{1-T^{\frac{1}{n_2}}} 
\end{eqnarray}
where $n_1$, $n_2$, $w_1$, $w_2$ are any positive integers which satisfy  
$\Delta:=w_2-w_1n_1n_2>0$.
The formula \eqref{5.4} of the quasihomogeneous case gives us the term
\begin{eqnarray}\label{6.5}
\sum_{i=0}^{k} {2k \choose 2i} B_{2i} B_{2(k-i)} 
\left( \frac{1}{(w_1n_1n_2)^{2i
-1}} - \frac{1}{w_2^{2i-1}} \right) \left( 1- \frac{1}{n_2^{2(k-i)-1}} \right).
\end{eqnarray}
We remark that we can factorise by $\Delta$ and $n_2-1$ and we get
\begin{eqnarray}\nonumber
(-B_0)B_{2k} \left( \frac{\sum_{j=0}^{2(k-1)}n_2^j}{n_2^{2k-1}} 
+ \frac{\sum_{j=0}^{2(k-1)}
w_2^j(w_1n_1n_2)^{2(k-1)-j}}{(w_1w_2n_1n_2)^{2k-1}} \right) \\
+
\sum_{i=1}^{k-1} {2k \choose 2i} B_{2i} B_{2(k-i)} 
\frac{\sum_{j=0}^{2(i-1)}w_2^
j(w_1n_1n_2)^{2(i-1)-j}}{(w_1w_2n_1n_2)^{2i-1}} 
\frac{\sum_{j=0}^{2(k-i-1)}n_2^j
}{n_2^{2(k-i)-1}}.\label{6.6}
\end{eqnarray}
This permits us to conclude.
\hfill $\qed$

\begin{remark}\label{t6.2}
For curves we expect in general to get
\begin{eqnarray}\label{6.7}
\Gamma^{Ber}(V^{sing}(f),2) =\Gamma^0 + \sum_{e \in Ed} \Gamma^e \Delta_e,
\end{eqnarray}
where $Ed$ is the set of edges of the Eisenbud and Neumann diagram of $f$ and 
$\Delta_e$ the determinant of the edge $e$. Because of the local situation, 
$\Delta_e$ is always positive.
In \cite{Br3} as well as above and in other cases, we have formulas as
\eqref{6.7} with $\Gamma^0$ of quasihomogeneous type.
The difficulty in proving the conjecture \ref{t1.2} (W) for other curve
singularities lies in the complexity of the coefficients $\Gamma^e$.
\end{remark}

\section{Compact complex manifolds}\label{c7}
\setcounter{equation}{0}

\noindent
The proof of theorem \ref{t1.6} will consist of three parts.
In the first part (A) we will derive the formula \eqref{7.4}
for $V^{mfd}(X)$. Motivated by it, we will define and discuss
the polynomials $q_{kj}(\nu,y_1,...,y_j)$ in part (B).
In part (C) we will prove theorem \ref{t1.6} and the formulas:
\begin{eqnarray}\label{7.1}
\Gamma_{2k}^{Ber}(V^{mfd}(X),\nu) &=& 
\sum_{j=0}^{\min(2k-1,n)} \int_X q_{kj}(n-\nu,c_1,...,c_j)\cdot c_{n-j}
\hspace*{0.5cm}\\
\hbox{ if }k\geq 1,&&\nonumber\\
\Gamma_{0}^{Ber}(V^{mfd}(X),\nu)&=&\int_Xc_n\label{7.2}
\end{eqnarray}
for any $\nu\in\C$.

After the proof we will make some remarks and finish with three examples.

\bigskip
{\bf (A)} Let $X$ be a compact complex manifold of dimension $n$.
Let $\alpha_1,...,\alpha_n$ be the Chern roots of the Chern classes
$c_1,...,c_n$, that is, $1+c_1+...+c_n=\prod_{j=1}^n(1+\alpha_j)$.
The Hirzebruch-Riemann-Roch theorem \cite{Hi}\cite{AS} gives
\begin{eqnarray}\label{7.3}
\chi(\Omega^p)&=&\int_X \left[Td(TM)\cdot ch(\Omega^p)\right]\\
&=& \int_X \left[ \left(\prod_{j=1}^n\frac{\alpha_j}{1-e^{-\alpha_j}}\right)
\cdot \sum_{j_1<...<j_p}e^{-\alpha_{j_1}-...-\alpha_{j_p}}\right].\nonumber
\end{eqnarray}
Here $Td(TM)$ is the Todd class of the tangent bundle $TM$ and 
$ch(\Omega^p)$ is the exponential Chern character of $\Omega^p$.
With \eqref{7.3} and \eqref{2.30} we can calculate the
generating function $V^{mfd}(X)$ of the higher moments $V^{mfd}_{2k}(X)$
which were defined in \eqref{1.6}.
\begin{eqnarray}\label{7.4}
&&V^{mfd}(X)\\ 
&=& \sum_{p=0}^n \chi_p e^{t(p-\frac{n}{2})} = 
e^{-\frac{n}{2}t}\sum_{p=0}^n \chi(\Omega^p)(-e^t)^p\nonumber \\
&=& e^{-\frac{n}{2}t} \int_X \left[
\left(\prod_{j=1}^n\frac{\alpha_j}{1-e^{-\alpha_j}}\right)
\cdot\prod_{j=1}^n (1-e^te^{-\alpha_j})\right] \nonumber\\
&=& \int_X\left[ \prod_{j=1}^n \left(\alpha_j
\frac{\sinh((\alpha_j-t)/2)}{\sinh(\alpha_j/2)}\right)\right]\nonumber\\
&=& \int_X\left[ \exp\left(\sum_{j=1}^n\left(
\Theta^{Ber}(\alpha_j) - \Theta^{Ber}(\alpha_j-t)\right)\right)\cdot
\prod_{j=1}^n(\alpha_j-t)\right].\nonumber
\end{eqnarray}

\bigskip
{\bf (B)}
Let $m\in \N_{\geq 1}$ be fixed. We will construct polynomials
$a_{k,2k-j}^{(m)},\ b_{kl}^{(m)}\in \Q[y_1,...,y_m]$
and $c_{kl}^{(m)},\ d_{kj}^{(m)}\in \Q[\nu,y_1,...,y_m]$.
They will all be quasihomogeneous of some degree (the second lower
index, $2k-j,l,l,j$) with respect to $y_1,...,y_m$, where 
$\deg y_j=j$. Those polynomials with weighted degree $\leq m$
will be independent of the choice of $m$; that means for example
in the case of $a_{k,2k-j}^{(m)}$ that $a_{k,2k-j}^{(m)}=a_{k,2k-j}^{(m')}$ 
for any $m'\geq 2k-j$.
At the end we will define $q_{kj}:=d_{2k,j}^{(j)}$.

Let $\sigma_{j}=\sigma_{j}(x_1,...,x_m)$, $j=1,...,m$, be the elementary
symmetric polynomials in $x_1,...,x_m$. 
For $k\geq 1$ and $1\leq j\leq 2k-1$ a unique polynomial
$a_{k,2k-j}^{(m)}\in \Q[y_1,...,y_m]$ exists such that
\begin{eqnarray}\label{7.5}
\sum_{i=1}^m(x_i^{2k}-(x_i-t)^{2k}+t^{2k}) 
= \sum_{j=1}^{2k-1} t^j\cdot a_{k,2k-j}^{(m)}(\sigma_1,...,\sigma_m).
\end{eqnarray}
It is quasihomogeneous of degree $2k-j$ with respect to $y_1,...,y_m$.
It is independent of $m$ in the sense described above if $2k-j\leq m$.

Because of \eqref{3.9}, for $k\geq 1$ and $l\geq 1$ unique polynomials
$b_{kl}^{(m)}\in \Q[y]$ exist which are quasihomogeneous of weighted
degree $l$ with respect to $y_1,...,y_m$ and which satisfy
\begin{eqnarray}\nonumber
&&\exp\left(\sum_{i=1}^m \left[ \Theta^{Ber}(x_i)-\Theta^{Ber}(x_i-t)
+\Theta^{Ber}(t)\right]\right)\\
&=& \exp\left(\sum_{k=1}^\infty \frac{-1}{2k}B_{2k}\frac{1}{(2k)!}
\sum_{j=1}^{2k-1}t^j\cdot a_{k,2k-j}^{(m)}(\sigma_1,...,\sigma_m)\right)
\nonumber\\
&=& 1 + \sum_{k=1}^\infty t^k\cdot  
\sum_{l=1}^\infty b_{kl}^{(m)}(\sigma_1,...,\sigma_m).\label{7.6}
\end{eqnarray}
The polynomials $b_{kl}^{(m)}$ with $l\leq m$ are independent of $m$.

For $k\in\N$ and $l\in \N$ unique polynomials 
$c_{kl}^{(m)}\in \Q[\nu,y_1,...,y_m]$ exist which are quasihomogeneous
of degree $l$ with respect to $y_1,...,y_m$ and which satisfy
\begin{eqnarray}\nonumber
&&\exp\left(\sum_{i=1}^m \left[ \Theta^{Ber}(x_i)-\Theta^{Ber}(x_i-t)
+\Theta^{Ber}(t)\right]\right)\exp(-\nu\Theta^{Ber}(t))\\
&&= \sum_{k=0}^\infty t^k\cdot \sum_{l=0}^\infty 
c_{kl}^{(m)}(\nu,\sigma_1,...,\sigma_m).\label{7.7}
\end{eqnarray}
The polynomials $c_{kl}^{(m)}$ with $l\leq m$ are independent of $m$.
Using \eqref{3.2} and \eqref{7.6} one sees
\begin{eqnarray}\label{7.8}
c_{k0}^{(m)}(\nu,y)&=& \frac{1}{k!}A_k(0,-\nu)
\ (=0 \hbox{ if }k\hbox{ is odd}),\\
c_{0l}^{(m)}(\nu,y)&=& 0 \quad \hbox{ if }l\geq 1,\label{7.9}\\
c_{kl}^{(m)}(\nu,y)&=& \sum_{j=0}^{k-1}\frac{1}{j!}A_j(0,-\nu)
\cdot b_{k-j,l}^{(m)}(y)\label{7.10}
\end{eqnarray}
if $k\geq 1$ and $l\geq 1$.
This implies especially
\begin{eqnarray}\label{7.11}
\deg_\nu c_{2k,0}^{(m)}=k,\quad
\deg_\nu c_{kl}^{(m)}\leq \left[\frac{k-1}{2}\right]\hbox{ if }
l\geq 1.
\end{eqnarray}

We define for $k\geq 1$ and $0\leq j\leq k-1$
\begin{eqnarray}\label{7.12}
d_{kj}^{(m)}(\nu,y):= k!(-1)^jc_{k-j,j}^{(m)}.
\end{eqnarray}
A simple calculation shows that the part of quasihomogeneous degree
$m$ with respect to $y_1,...,y_m$ in 
\begin{eqnarray}\label{7.13}
\left(\sum_{k=0}^\infty t^k \cdot\sum_{l=0}^\infty 
c_{kl}^{(m)}(\nu,y)\right)\cdot \left(\sum_{i=0}^my_{m-i}(-t)^i\right)
\end{eqnarray}
(with $y_0:=1$) is
\begin{eqnarray}\label{7.14}
y_m+\sum_{k=1}^\infty \frac{1}{k!}t^k\cdot \sum_{j=0}^{\min(k-1,m)} 
y_{m-j}d_{kj}^{(m)}(\nu,y).
\end{eqnarray}
Finally, we define for $k\geq 1$ and $1\leq j\leq 2k-1$ 
the polynomials $q_{kj}(\nu,y)$ by
\begin{eqnarray}\label{7.15}
q_{kj}:=d_{2k,j}^{(j)}.
\end{eqnarray}
They are quasihomogeneous of degree $j$ with respect to $y_1,...,y_m$;
the degrees with respect to $\nu$ satisfy because of \eqref{7.11}
\begin{eqnarray}\label{7.16}
\deg_\nu q_{k0}=k,\quad \deg_\nu q_{kj}\leq k-1-\left[\frac{j}{2}\right]
\hbox{ if }j\geq 1.
\end{eqnarray}

\bigskip
{\bf (C)} In \eqref{7.4} the last factor is (with $c_0:=1$)
\begin{eqnarray}\label{7.17}
\prod_{j=1}^n(\alpha_j-t)= \sum_{i=0}^n c_{n-i}(-t)^i.
\end{eqnarray}
$\Gamma^{Ber}(V^{mfd}(X),n-\nu)$ contains only even powers of $t$.
Combining \eqref{7.4}, \eqref{7.7}, \eqref{7.13} and \eqref{7.14},
we find
\begin{eqnarray}\label{7.18}
&&\Gamma^{Ber}(V^{mfd}(X),n-\nu) \\
&=& \int_Xc_n+ \sum_{k=1}^\infty
\frac{1}{(2k)!}t^{2k} \sum_{j=0}^{\min(2k-1,m)} 
\int_X c_{m-j}q_{kj}(\nu,c_1,...,c_j)\nonumber.
\end{eqnarray}
This shows \eqref{7.1}, \eqref{7.2} and theorem \ref{t1.6}.

\begin{remarks}
a) The formula \eqref{1.33} for  $V^{mfd}_2(X)$ was calculated in \cite{LW}
and \cite{Bo}. Calculations with some ressemblance to those
in (A) can be found in \cite[\S 3]{Sal}.

b) By \eqref{7.2}, for any compact complex manifold $X$
$\Gamma_0^{Ber}(V^{mfd}(X),n)=\int_X c_n.$
On the other hand, analogously to \eqref{1.22} and \eqref{1.23},
the sequence of numbers
\begin{eqnarray}\label{7.19}
(-1)^k \Gamma^{Ber}_{2k}(V^{mfd}(X),n)\cdot 
\frac{(2\pi)^{2k}\cdot \Gamma(n)} {2\cdot (2k)!\cdot (2k)^{n-1}}
\end{eqnarray}
tends with $k\to \infty$ to $(-1)^n\sum_p\chi_p=(-1)^n\int_X c_n$.
Therefore for odd $n$ and $\int_Xc_n\neq 0$ the analogue of the
conjectures \ref{t1.2} is never satisfied.

c) The example $X=\P^n$ below shows a different rule for the signs
if $2k< n$. The example of a K3 surface below shows a behaviour
analogous to quasihomogeneous singularities.
One could try to classify the compact complex manifolds
according to the behaviour of the signs of 
$(-1)^k \Gamma^{Ber}_{2k}(V^{mfd}(X),n)$.
\end{remarks}

\begin{examples}\label{t7.2}
a) $X=\P^n$: We use N\"orlunds notation 
$B_k^{(\nu)}(x)$ of the generalized Bernoulli polynomials, see \eqref{3.3},
because the generalized Bernoulli numbers $B_k^{(\nu)}(0)$ ($k,\nu\in\N$)
will play a role.
\begin{eqnarray}\label{7.20}
V^{mfd}(\P^n)= e^{-\frac{n}{2}t}\sum_{p=0}^n e^{tp} 
= e^{-\frac{n}{2}t}\cdot \frac{e^{t(n+1)}-1}{e^t-1}
\end{eqnarray}
and 
\begin{eqnarray}\label{7.21}
&& \Gamma^{Ber}(V^{mfd}(\P^n),n)= V^{mfd}(\P^n)\cdot 
\left(\frac{te^{\frac{1}{2}t}}{e^t-1}\right)^n \\
&=& \frac{1}{t}\left[ \frac{t^{n+1}e^{t(n+1)}}{(e^t-1)^{n+1}}
-\frac{t^{n+1}}{(e^t-1)^{n+1}}\right]\nonumber\\
&=& \frac{1}{t} \sum_{k=0}^\infty \frac{1}{k!}t^k\left[ B_k^{(n+1)}(n+1)
-B_k^{(n+1)}(0)\right]\nonumber\\
&=& \sum_{k=0}^\infty \frac{1}{(2k)!}t^{2k}\cdot \frac{-2}{2k+1}
B^{(n+1)}_{2k+1}(0).\nonumber 
\end{eqnarray}
We used the (anti-)symmetry \eqref{3.17}. Now \eqref{3.22} and a special
case of \eqref{3.16} show (see also \cite[p 148]{No2})
\begin{eqnarray}\label{7.22}
(x-1)...(x-n)=B_n^{(n+1)}(x)=\sum_{s=0}^n {n\choose s}x^sB_{n-s}^{(n+1)}(0).
\end{eqnarray}
Therefore, for $s<n+1$ 
the sign of the generalized Bernoulli number $B^{(n+1)}_{s}(0)$ is $(-1)^s$
Thus, $(-1)^k\Gamma^{Ber}_{2k}(V^{mfd}(\P^n),n)= (-1)^k\frac{-2}{2k+1}
B^{(n+1)}_{2k+1}(0)$ has for $2k<n$ the sign $(-1)^k$.
This behaviour is completely different from that for large $k$ and that
in the conjectures \ref{t1.2}.

\bigskip
b) $X$ a K3 surface: 
\begin{eqnarray}\label{7.23}
V^{mfd}(K3)&=& e^{-t}(2+20e^t+2e^{2t})\\
&=& 2\cdot V^{mfd}(\P^2) + 18,\nonumber
\end{eqnarray}
\begin{multline}
\Gamma^{Ber}(V^{mfd}(K3),2)= 2\Gamma^{Ber}(V^{mfd}(\P^2),2)+ 18
\cdot\left(\frac{te^{\frac{1}{2}t}}{e^t-1}\right)^2\label{7.24}\\
= \sum_{k=0}^\infty \frac{1}{(2k)!}t^{2k}
\left[\frac{-4}{2k+1}B_{2k+1}^{(3)}(0) + 18\cdot B_{2k}^{(2)}(1)\right].
\end{multline}
One can calculate \eqref{7.25} and \eqref{7.26} with \eqref{3.21},
and then \eqref{7.27} with \eqref{3.20};
\begin{eqnarray}\label{7.25}
B_k^{(2)}(0)&=& (1-k)B_k-kB_{k-1},\\
B_k^{(3)}(0)&=& (2-k)B_k^{(2)}(0)-2kB_{k-1}^{(2)}(0)\nonumber\\
&=& \frac{1}{2}(k-2)(k-1)B_k+ \frac{3}{2}(k-2)kB_{k-1}
\label{7.26}\\
&&+(k-1)kB_{k-2},\nonumber\\
B_{k}^{(2)}(1)&=&B_{k}^{(2)}(0)+kB_{k-1} = (1-k)B_{k}.\label{7.27}
\end{eqnarray}
Therefore 
\begin{eqnarray}\nonumber
&&(-1)^k\Gamma^{Ber}_{2k}(V^{mfd}(K3),2) \\
&=& 24(2k-1)(-1)^{k-1}B_{2k}
+ (-1)^{k-1}8kB_{2k-1}\nonumber \\
&=&\left\{\begin{array}{rrr} 0&\hbox{ if }&k=1,\\
  24(2k-1)(-1)^{k-1}B_{2k} >0 &\hbox{ if }&k\neq 1.\end{array}\right.
\label{7.28}
\end{eqnarray}
This behaviour is similar to that of a quasihomogeneous singularity.

\bigskip
c) $X$ a Riemann surface of genus $g$:
\begin{eqnarray}\label{7.29}
V^{mfd}(X)&=&(1-g)V^{mfd}(\P^1),\\
\Gamma^{Ber}(V^{mfd}(X),1) &=&(1-g)\Gamma^{Ber}(V^{mfd}(\P^1),1)\label{7.30}\\
&=& (1-g)\sum_{k=0}^\infty \frac{1}{(2k)!}t^{2k}\cdot 2B_{2k}.\nonumber
\end{eqnarray}
We used \eqref{7.25}. For $g\geq 2$ $(-1)^k\Gamma^{Ber}_{2k}(V^{mfd}(X),1)$
is positive if $k\geq 1$, but negative if $k=0$.
\end{examples}

\end{document}